\documentclass[twoside,11pt]{amsart}

\usepackage[a4paper,hmargin={2.20cm,2.20cm},vmargin={2cm,2cm},heightrounded,marginparwidth=2.1cm, marginparsep=0.05cm]{geometry}
\usepackage[pdftex]{hyperref}
\usepackage{bbm}
\usepackage{xcolor}
\usepackage[all,arc,cmtip]{xy}
\usepackage{amstext}
\usepackage{amsfonts}
\usepackage{amsmath}
\usepackage{amssymb}
\usepackage{amsthm}
\usepackage{mathtools}
\usepackage{rotating}
\usepackage{cmll}
\usepackage{stmaryrd}
\usepackage{pigpen}
\usepackage{slashbox}
\usepackage{mathcomp}
\usepackage{tikz-cd}
\usepackage{lineno}
\usepackage{scalerel}

\def\DS{\displaystyle}
\newcommand{\categ}[1]{\mathsf{#1}}
\newcommand{\Cat}[1]{#1\text{-}\categ{Cat}}
\newcommand{\Rell}[1]{#1\text{-}\categ{Rel}}
\def\rel{{\longrightarrow\hspace*{-2.8ex}{\mapstochar}\hspace*{2.8ex}}}
\def\hom{{\rm hom}}

\usepackage{chngcntr}
\counterwithin*{equation}{section}

\def\Set{\categ{Set}}

\def\Top{\categ{Top}}

\def\App{\categ{App}}
\def\Ord{\categ{Ord}}
\def\Met{\categ{Met}}
\def\UltMet{\categ{UltMet}}  

\def\NA{\categ{NA}}  
\def\Qs{\categ{Qs}}
\def\cX{\categ{X}}
\def\one{\categ{1}}
\def\A{\categ{A}}

\def\C{\mathcal{C}}
\def\mD{\mathcal{D}}

\def\x{\mathfrak{x}}
\def\y{\mathfrak{y}}
\def\w{\mathfrak{w}}
\def\z{\mathfrak{z}}
\def\xx{\mathfrak{X}}
\def\p{\mathfrak{p}}
\def\q{\mathfrak{q}}
\def\v{\mathfrak{v}}

\def\U{\mathbbm{U}}
\def\I{\mathbbm{I}}
\def\T{\mathbbm{T}}

\newcommand{\V}{\categ{V}}
\newcommand{\PP}{\categ{P}}
\def\two{\categ{2}}

\title{Compactly generated spaces and quasi-spaces in topology}
\author{Willian Ribeiro}
\address{CMUC, Department of Mathematics, University of Coimbra, 3001-501 Coimbra, Portugal}
\email{willian.ribeiro.vs@gmail.com}

\thanks{Research supported by Centro de Matem\'{a}tica da Universidade de Coimbra {\bf --} UID/MAT/00324/2019 and by the FCT PhD Grant PD/BD/128059/2016, funded by the Portuguese Government through FCT/MCTES and co-funded by the European Regional Development Fund through the Partnership Agreement PT2020.}

\keywords{$(\T,\V)$-categories, compact and Hausdorff space, compactly generated space, Alexandroff space, cartesian closedness, quasi-space}

\subjclass[2010]{18B30, 18D15, 54A05, 54B30, 54C35, 54D50}

\date{\today}

\begin{document}

\begin{abstract}
The notions of compactness and Hausdorff separation for generalized enriched categories allow us, as classically done for the category $\Top$ of topological spaces and continuous functions, to study {\it compactly generated spaces} and {\it quasi-spaces} in this setting. Moreover, for a class $\C$ of objects we generalize the notion of {\it $\C$-generated spaces}, from which we derive, for instance, a general concept of {\it Alexandroff spaces}. Furthermore, as done for $\Top$, we also study, in our level of generality, the relationship between compactly generated spaces and quasi-spaces. 
\end{abstract}

\maketitle

\tableofcontents


\section*{Introduction}

Amongst the classical approaches to the inconvenience of non-cartesian closedness of $\Top$ {\bf --} restricting ourselves to a subcategory thereof or including it in a supercategory {\bf --} in this paper we turn our attention to the subcategory of {\it compactly generated spaces} and the supercategory of {\it quasi-topological spaces}. 

Compactly generated spaces were widely studied in the second half of the past century, having as early references the works of Kelley \cite{MR0070144} and Gale \cite{MR0036503}, who indicates Hurewicz as the first to define the notion. The main motivation for its investigation was the search for a category of spaces convenient for homotopy theory and algebraic topology \cite{MR0210075, MR0300277}. The full subcategory of $\Top$ of compactly generated spaces, that we will denote by $\Top_{_{\C}}$, where $\C$ stands for the class of compact Hausdorff topological spaces, is well known to be cartesian closed and coreflective in $\Top$; it is actually the {\it coreflective hull} of $\C$ in $\Top$, that is, the smallest coreflective subcategory of $\Top$ that contains $\C$, from where we conclude that compactly generated topological spaces coincide with colimits of compact Haudorff topological spaces. This coreflection is also shown in \cite{Mac71}. In order to establish our results, we directly follow the approach of Escard\'o, Lawson and Simpson \cite{MR2080286}, which is centered in the concept of a generating class of objects, in this case, the compact Hausdorff ones.

Concerning quasi-topological spaces, as indicated by Spanier \cite{MR0144300}, they were introduced with the intention of constructing internal function spaces in the resulting category $\Qs\Top$ of quasi-topological spaces and quasi-continuous maps. The latter supercategory of $\Top$, that is convenient for homotopy theory, was suitably used in the results of Booth \cite{MR0315710} and Day \cite{MR642246}. Although carrying a size {\it illegitimacy} demonstrated by Herrlich and Rajagopalan \cite{MR708599}, the category of quasi-topological spaces provided a useful tool in current works as the ones by Dadarlat and Meyer \cite{MR2920847} and Browne \cite{ETSGCAB} on {\it $E$-theory of $C^{*}$-algebras}, which share a common topic with a past work of Dubuc and Porta \cite{MR566280}; we also refer to Petrakis' PhD thesis \cite{ctbsp} that highlights the relation between quasi-topological spaces and Bishop spaces. Moreover, the category $\Qs\Top$ also serves as a paradigm for the results of Escard\'o and Xu \cite{MR3107257}, and Dubuc and Espa{\~n}ol \cite{MR555548, Dubuc06quasitopoiover} who presented a much more general context of quasi-topologies using the notion of {\it Grothendieck topologies}.    

The main goal of the present work is to carry these two concepts {\bf --} compactly generated spaces and quasi-spaces {\bf --} from $\Top$ to $\Cat{(\T,\V)}$. The $(\T,\V)${\it -setting} \cite{MR1957813, MR1990036} (see also \cite{MR3307673}) has been broadly studied and applied to problems in topology \cite{MR2355608, MR2107395, MR2729224, MR3987969}, for it provides a unified effective way of studying classical categories from Analysis and Topology such as categories of ordered, metric, topological, and approach spaces.

This setting also allows for a generalization of compact and Hausdorff spaces, which is fundamental for our discussion. We recall this generalization in the first section, as well as some essential facts about topological functors, and necessary background of $(\T,\V)$-spaces and $(\T,\V)$-continuous maps. In Section \ref{sec2} we develop the concept of $\C$-generated $(\T,\V)$-spaces generalizing the one established for $\Top$ in \cite{MR2080286}; in particular, we get the notions of {\it compactly generated $(\T,\V)$-spaces} and {\it Alexandroff $(\T,\V)$-spaces}. We finish with Section \ref{sec3} presenting the quasi-spaces in our general context and studying a relationship between the two main notions, extending the case of $\Top$ due to Day \cite{rsqsksd}, whose work we follow closely. Examples are presented throughout the paper.


\section{Prerequisites}\label{sec1}


\subsection{A comment on topological functors}

Concerning topological functors, we recall some facts that can be found in the literature. Let $\A$ and $\cX$ be categories and $|\text{-}|\colon\A\to\cX$ be a functor. A sink $(f_{_{i}}\colon A_{_{i}}\to A)_{_{i\in I}}$ of morphisms of $\A$ is {\it $|\text{-}|$-final} if, for every sink $(g_{_{i}}\colon A_{_{i}}\to B)_{_{i\in I}}$ in $\A$ and every morphism $s\colon|A|\to|B|$ of $\cX$ with $s\cdot|f_{_{i}}|=|g_{_{i}}|$, there exists a unique $t\colon A\to B$ in $\A$ with $|t|=s$ and, for all $i\in I$, $t\cdot f_{_{i}}=g_{_{i}}$.
\begin{align*}\xymatrix{|A_{_{i}}| \ar[r]^{|f_{_{i}}|} \ar[rd]_{|g_{_{i}}|} & |A| \ar[d]^{s} & A \ar@{..>}[d]^{\exists \ ! \ t} \\ & |B| & B}\end{align*}
The dual concept is that of a {\it $|\text{-}|$-initial source}. The functor $|\text{-}|$ is {\it topological} if every {\it $|\text{-}|$-structured sink}, that is, a sink of the form $(f_{_{i}}\colon |A_{_{i}}|\to X)_{_{i\in I}}$ in $\cX$, with $A_{_{i}}\in\A$, admits a {\it $|\text{-}|$-final lifting}, or, equivalently, if every {\it $|\text{-}|$-structured source} $(f_{_{i}}\colon X\to|A_{_{i}}|)_{_{i\in I}}$ in $\cX$ admits a {\it $|\text{-}|$-initial lifting} \cite[Theorem 21.9]{MR1051419}.

Observe that, as usual, we do not make any assumption of smallness on the entity $I$ of indexes. However, in our case, we will be dealing with a {\it fibre-small} (the fibre of each object in the codomain is a set) forgetful functor, hence we will employ the fact that each $|\text{-}|$-final lifting of a sink \linebreak $(g_{_{i}}\colon |A_{_{i}}|\to X)_{_{i\in I}}$ is actually the $|\text{-}|$-final lifting of a sink $(g_{_{j}}\colon |A_{_{j}}|\to X)_{_{j\in J}}$, with $J$ a set contained in $I$ \cite[Proposition 21.34]{MR1051419}. 


\subsection{\texorpdfstring{$(\T,\V)$}{(T,V)}-spaces and \texorpdfstring{$(\T,\V)$}{(T,V)}-continuous maps}\label{sub_tvcat}

We introduce briefly the $\Cat{(\T,\V)}$ setting for our purposes, and refer the reader to the references \cite{MR1990036, MR1957813, MR3307673} for details.

Let $\V=(\V,\otimes,k)$ be a unital commutative quantale (see, for instance, \cite[II-1.10]{MR3307673}), which is also a {\it Heyting algebra}, so that the operation infimum $\wedge$ has a right adjoint. Consider the order-enriched category $\Rell{\V}$: objects are sets and morphisms are $\V$-relations (or $\V$-matrices) $r\colon X\rel Y$, which are $\V$-valued maps $r\colon X\times Y\to\V$; the {\it relational composition} of $r\colon X\rel Y$ and $s\colon Y\rel Z$ is given by: for each $(x,z)\in X\times Z$, $s\cdot r(x,z)=\DS\bigvee_{_{y\in Y}}\left(r(x,y)\otimes s(y,z)\right)$, and the order between $\V$-relations is defined {\it componentwise}. There exists an {\it involution} given by transposition: for each $r\colon X\rel Y$, $r^{\circ}\colon Y\rel X$ is given by, for each $(y,x)\in Y\times X$, $r^{\circ}(y,x)=r(x,y)$. Denoting the bottom element of the complete lattice $\V$ by $\bot$, each map $f\colon X\to Y$ can be seen as a $\V$-relation $f\colon X\rel Y$: for $(x,y)\in X\times Y$,
\begin{align*}f(x,y)=\left\{\begin{array}{ll}k, & \text{if $f(x)=y$} \\ \bot, & \text{otherwise.}\end{array}\right.\end{align*}

Let $\T=(T,m,e)\colon\Set\to\Set$ be a monad satisfying the {\it Beck-Chevalley condition}, (BC) for short, that is, $T$ preserves weak pullbacks and the naturality squares of $m$ are weak pullbacks \cite{MR3175322}. We fix a lax extension of $\T$ to $\Rell{\V}$, again denoted by $\T$, so that $T\colon\Rell{\V}\to\Rell{\V}$ is a lax functor and the natural transformations $m$ and $e$ become oplax: for each $\V$-relation $r\colon X\rel Y$, 
\begin{align*}\xymatrix{X \ar[r]^{e_{X}} \ar[d]|(0.45){\object@{|}}_(0.45){r} \ar@{}[dr]|{\le} & TX \ar[d]|(0.45){\object@{|}}^(0.45){\hat{T}r} & T^{2}X \ar[l]_{m_{X}} \ar[d]|(0.45){\object@{|}}^(0.45){\hat{T}^{2}r} \ar@{}[dl]|{\ge} \\ Y \ar[r]_{e_{Y}} & TY. & T^{2}Y \ar[l]^{m_{Y}}}\end{align*}
We assume that this lax extension is {\it flat}, that is, for each set $X$, $T1_{_{X}}=1_{_{TX}}$, and, moreover, that it commutes with involution, i.e., for each $r\colon X\rel Y$ in $\Rell{V}$, $T(r^{\circ})=(Tr)^{\circ}$. Then we have a lax monad on $\Rell{\V}$ in the sense of \cite{MR2116322} and $\Cat{(\T,\V)}$ is defined as its category of Eilenberg-Moore lax algebras. Hence objects are pairs $(X,a)$, where $X$ is a set and $a\colon TX\rel X$ is a {\it reflexive} and {\it transitive} $\V$-relation, so that the diagram
\begin{align*}\xymatrix{X \ar[r]^{e_{X}} \ar@/_1pc/[dr]_{1_{X}} & TX \ar[d]|(0.45){\object@{|}}^(0.45){a} & T^{2}X \ar[l]|{\object@{|}}_{Ta} \ar[d]^{m_{X}} \ar@{}[dl]|{\le} \\ \ar@{}[ur]|{\le} & X & TX \ar[l]|(0.45){\object@{|}}^(0.45){a}}\end{align*}
is lax commutative; such pairs are called {\it $(\T,\V)$-categories} or {\it $(\T,\V)$-spaces}; and a morphism from $(X,a)$ to $(Y,b)$ is a map $f\colon X\to Y$ such that the square
\begin{align*}\xymatrix{TX \ar[d]|(0.45){\object@{|}}_(0.45){a} \ar[r]^{Tf} \ar@{}[dr]|{\le} & TY \ar[d]|(0.45){\object@{|}}^(0.45){b} \\ X \ar[r]_{f} & Y}\end{align*}
is lax commutative; such a map is called a {\it $(\T,\V)$-functor} or a {\it $(\T,\V)$-continuous} map. When this diagram is strictly commutative, what is equivalent to $a=f^{\circ}\cdot b\cdot Tf$, $f$ is said to be {\it fully faithful}. The $\V$-relation $a\colon TX\rel X$ of the $(\T,\V)$-space $(X,a)$ is referred usually as the {\it $(\T,\V)$-structure} of $X$. For a $(\T,\V)$-space $(X,a)$, each subset $A\subseteq X$ can be endowed with a {\it subspace $(\T,\V)$-structure} $a_{_{A}}=i_{_{A}}^{\circ}\cdot a\cdot Ti_{_{A}}\colon TA\rel A$, where $i_{_{A}}\colon A\hookrightarrow X$ is the inclusion map, which becomes a fully faithful map $i_{_{A}}\colon(A,a_{_{A}})\to(X,a)$. 

The forgetful functor $|$-$|\colon\Cat{(\T,\V)}\to\Set$ is topological \cite{MR1990036,MR1957813} and fibre-small: for each set $X$, a $(\T,\V)$-structure $a$ on $X$ is an element of $\Rell{\V}(TX,X)=\Set(TX\times X,\V)$. 

We also need to assume the fairly restrictive condition that each constant map between $(\T,\V)$-spaces is continuous, what implies, as we show below, that, in particular, the quantale $\V$ is {\it integral}, i.e., $k=\top$ is the top element of $\V$. \vspace{0.2cm}

\noindent{\bf Lemma} {\it The following statements are equivalent.

\noindent(i) Any constant map $y_{_{0}}\colon(X,a)\to(Y,b)$ between $(\T,\V)$-spaces is continuous.

\noindent(ii) If $(\one,c)$ is a $(\T,\V)$-space, where $\one$ denotes a singular set $\{*\}$, then, for each $\z\in T\one$, $c(\z,*)=\top$.

\noindent(iii) $k=\top$ and $T1=1$.} \vspace{0.2cm}

\noindent{\it Proof.} (i)$\Leftrightarrow$(ii) Let $c$ be a $(\T,\V)$-structure on $\one$. The identity map $1_{_{\one}}\colon(\one,\top)\to(\one,c)$ is constant, hence it is continuous, so $\top\leq c$. Conversely, for the constant map $y_{_{0}}\colon(X,a)\to(Y,b)$, consider the factorization
\begin{align*}\xymatrix{(X,a) \ar[rr]^{y_{_{0}}} & & (\one,b_{_{\one}}) \ar[rr]^{i_{_{\one}}} & & (Y,b),}\end{align*}
where $\one=\{y_{_{0}}\}\subseteq Y$ and $b_{_{\one}}$ is the subspace $(\T,\V)$-structure. By hypothesis, $b_{_{\one}}$ is constantly equal to $\top$, what implies that $y_{_{0}}\colon(X,a)\to(\one,b_{_{\one}})$ is continuous, whence the composite $y_{_{0}}\colon(X,a)\to(Y,b)$ is continuous. 

\noindent (ii)$\Leftrightarrow$(iii) Since $T$ is flat, the {\it discrete $(\T,\V)$-structure} on $\one$ is given by $e_{_{\one}}^{\circ}$, so we have the {\it discrete $(\T,\V)$-space} $(\one,e_{_{\one}}^{\circ})$ \cite[III-3.2]{MR3307673}. Assuming (ii), for each $\z\in T\one$, $e_{_{\one}}^{\circ}(\z,*)=\top$. In particular, for $\z=e_{_{\one}}(*)$,
\begin{gather*}k=e_{_{\one}}^{\circ}(e_{_{\one}}(*),*)=\top.
\shortintertext{Then, for each $\z\in T\one$,}
e_{_{\one}}^{\circ}(\z,*)=\top=k \ \Longleftrightarrow \ \z=e_{_{\one}}(*),\end{gather*}
whence $T\one=\one$. On the other hand, for $\V$ integral and $T1=1$, one readily checks condition (ii). \qed \vspace{0.2cm}

Under those conditions, $\Cat{(\T,\V)}$ is a topological category in the classical sense of \cite{MR0460414}, that is, there exist $|\text{-}|$-initial $(\T,\V)$-structures, with the forgetful functor $|\text{-}|\colon\Cat{(\T,\V)}\to\Set$ being fibre-small, and on a singleton set there exists exactly one $(\T,\V)$-structure. To give examples we will consider $\V$ as the following integral and {\it completely distributive} quantales (see, for instance, \cite[II-1.11]{MR3307673}): $\two=(\{\bot,\top\},\wedge,\top)$, $\categ{P}_{_{+}}=([0,\infty]^{\rm{op}},+,0)$, $\categ{P}_{_{{\rm max}}}=([0,\infty]^{\rm{op}},{\rm max},0)$, and \linebreak $[0,1]_{_{\odot}}=([0,1],\odot,1)$, where $\odot$ is the {\it \L{}ukasiewicz tensor} given by $u\odot v=\max(0,u+v-1)$, for each $u,v\in[0,1]$; in addition to $\T$ as: \\
$\bullet$ the identity monad $\mathbb{I}=({\rm Id},\one,\one)$ on $\Set$ extended to the identity monad on $\Rell{\V}$, and \\ 
$\bullet$ the ultrafilter monad $\U$ with the Barr extension to $\Rell{\V}$ \cite[IV-2.4.5]{MR3307673}.

In the categories of the following table, all constant maps are $(\T,\V)$-continuous.
\begin{align}\label{tab1}\text{\begin{tabular}{|c|c|c|c|c|} 
\hline
\backslashbox{$\T$}{$\V$} & $\two$ & $\categ{P}_{_{+}}$ & $\categ{P}_{_{{\rm max}}}$ & $[0,1]_{_{\odot}}$ \\ \hline
$\I$                      & $\Ord$ & $\Met$             & $\UltMet$                  & $\categ{B}_{_{1}}\Met$ \\ \hline
$\U$                      & $\Top$ & $\App$             & $\NA\text{-}\App$          & $\Cat{(\U,[0,1]_{_{\odot}})}$ \\ 
\hline\end{tabular}}\end{align}
$\bullet$ $\Ord$ is the category of pre-ordered spaces and monotone maps; \\
$\bullet$ $\Met$ is the category of Lawvere's generalized metric spaces and non-expansive maps \cite{MR1925933}; \\
$\bullet$ $\UltMet$ is the full subcategory of $\Met$ of ultrametric spaces \cite[III-Exercise 2.B]{MR3307673}; \\
$\bullet$ $\categ{B}_{_{1}}\Met$ is the category of {\it bounded-by-1} metric spaces and non-expansive maps (see \cite{MR3987969}); \\
$\bullet$ $\Top$ is the usual category of topological spaces and continuous functions; \\
$\bullet$ $\App$ is the category of Lowen's approach spaces and contractive maps \cite{MR1472024}; \\
$\bullet$ $\NA\text{-}\App$ is the full subcategory of $\App$ of non-Archimedean approach spaces \cite{MR3731477}. 


\subsection{Injective and exponentiable \texorpdfstring{$(\T,\V)$}{(T,V)}-spaces}\label{sub_inj_exp} 

Details on the following concepts can be found in \cite{MR1957813, MR2491799, MR2729224}.  

For each $(X,a),(Y,b)$ in $\Cat{(\T,\V)}$, there exists an induced (pre-)order on the set \linebreak $\Cat{(\T,\V)}((X,a),(Y,b))$ of $(\T,\V)$-continuous maps from $(X,a)$ to $(Y,b)$ given by, for each \linebreak $f,g\colon(X,a)\to(Y,b)$, 
\begin{align}\label{order_maps}f\leq g \ \Longleftrightarrow \ \forall x\in X, \ k\leq b(e_{_{Y}}(f(x)),g(x)).\end{align}
Denoting by $f\simeq g$ when $f\leq g$ and $g\leq f$, the $(\T,\V)$-space $(Y,b)$ is {\it separated} if, for every $(\T,\V)$-space $(X,a)$ and every $f,g\colon(X,a)\to(Y,b)$, $f\simeq g$ implies $f=g$, that is, for every $(X,a)$, the order (\ref{order_maps}) is anti-symmetric. Furthermore, the space $(Y,b)$ is separated precisely when the following order on $Y$ is anti-symmetric: for each $y,y'\in Y$,
\begin{align}\label{order_space}y\leq y' \ \Longleftrightarrow \ k\leq b(e_{_{Y}}(y),y'),\end{align}
i.e., when $y\leq y'$ in the order (\ref{order_maps}) for $y,y'\colon\one\to Y$ the morphisms induced by $y$ and $y'$, respectively. The full subcategory of $\Cat{(\T,\V)}$ of separated spaces, which is denoted by $\Cat{(\T,\V)}_{_{\rm sep}}$, is closed under mono-sources. 

A space $(Z,c)$ is {\it injective} if, for each fully faithful map $y\colon(X,a)\to(Y,b)$ and $(\T,\V)$-continuous map $f\colon(X,a)\to(Z,c)$, there exists a $(\T,\V)$-continuous map $\hat{f}\colon(Y,b)\to(Z,c)$, called an {\it extension of $f$ along $y$}, such that $\hat{f}\cdot y\simeq f$.
\begin{align}\label{diag_ext}\xymatrix{X \ar[rr]^{y} \ar[rd]_{f} \ar@{}[drr]^(0.45){\simeq} & & Y \ar[ld]^{\hat{f}} \\ & Z & }\end{align}
Observe that when considering separated $(\T,\V)$-spaces, injectivity assumes its usual notion. In \cite{MR2729224}, injective $(\T,\V)$-spaces are characterized as the ones satisfying a {\it cocompleteness} condition. 

As usual, a $(\T,\V)$-space $(X,a)$ is said to be {\it exponentiable} in $\Cat{(\T,\V)}$ if the functor \linebreak $\text{--}\times(X,a)\colon\Cat{(\T,\V)}\to\Cat{(\T,\V)}$ has a right adjoint. In order to recall conditions under which injective $(\T,\V)$-spaces are exponentiable, which are established in \cite{MR3987969}, we restrict ourselves to the case when the extension of $T$ to $\Rell{\V}$ is fully determined by a $\T$-algebra structure $\xi\colon T\V\to\V$, so we are in the setting of {\it strict topological theories} \cite{MR2355608}; such extensions are characterized in \cite{MR3330902} as the {\it algebraic extensions}. In particular, the following diagrams are commutative:
\begin{align}\label{diag_top_the}\xymatrix{\V \ar[r]^{e_{_{\V}}} \ar[dr]_{1_{_{\V}}} & T\V \ar[d]^(0.45){\xi} & T^{2}\V \ar[l]_{T\xi} \ar[d]^{m_{_{\V}}} \\ & \V & T\V \ar[l]^(0.45){\xi}}\qquad\xymatrix{T(\V\times\V) \ar[rr]^{T(\otimes)} \ar[d]_{\left\langle\xi\cdot T\pi_{_{1}},\xi\cdot T\pi_{_{2}}\right\rangle} & & T\V \ar[d]^{\xi} & T\one \ar[l]_{T(k)} \ar[d]^{!} \\ \V\times\V \ar[rr]_{\otimes} & & \V & \one, \ar[l]^{k}}\end{align}
where $\pi_{_{1}}$ and $\pi_{_{2}}$ are the product projections $\V\times\V\to\V$, and $\left\langle\xi\cdot T\pi_{_{1}},\xi\cdot T\pi_{_{2}}\right\rangle$ is the unique map rendering the diagram below commutative.
\begin{align*}\xymatrix@R=1.5em{ & & T(\V\times\V) \ar[ddll]_{\xi\cdot T\pi_{_{1}}} \ar[dd]|{\left\langle \xi\cdot T\pi_{_{1}},\xi\cdot T\pi_{_{2}}\right\rangle} \ar[ddrr]^{\xi\cdot T\pi_{_{2}}} & & \\ & & & & \\ \V & & \V\times\V \ar[ll]^{\pi_{_{1}}} \ar[rr]_{\pi_{_{2}}} & & \V}\end{align*} 
The extension of $T$ to $\Rell{\V}$ is given by, for each $r\colon X\rel Y$, $\x\in TX$, $\y\in TY$,
\begin{align}\label{extension}Tr(\x,\y)=\bigvee\{\xi\cdot T\vec{r}(\w) \ | \ \w\in T(X\times Y), \ T\pi_{_{X}}(\w)=\x, \ T\pi_{_{Y}}(\w)=\y\},\end{align}
where $\pi_{_{X}}$ and $\pi_{_{Y}}$ are the product projections from $X\times Y$ to $X$ and $Y$, respectively \cite[Definition 3.4]{MR2355608}, and we adopt the notation from \cite{MR3330902}: $T\vec{r}\colon T(X\times Y)\to T\V$ is the image of the map $\vec{r}\colon X\times Y\to\V$ by the functor $T$, while $Tr$ continues to have the meaning so far used, that is, it is the image of the $\V$-relation $r\colon X\rel Y$ by the extension of the functor $T$ to $\Rell{\V}$. 

In this context, $\V$ can be endowed with a $(\T,\V)$-structure $\hom_{_{\xi}}\colon T\V\rel\V$ given by the composite
\begin{align*}\xymatrix{T\V \ar[rr]^{\xi} & & \V \ar[rr]|{\object@{|}}^{\hom} & & \V,}\end{align*}
where $\hom\colon\V\times\V\to\V$ is the left adjoint of the tensor operator: for each $u,v,w\in\V$, 
\begin{align*}u\otimes v\leq w \ \Longleftrightarrow \ u\leq\hom(v,w),\end{align*}
and it is given by, for each $u,v\in\V$, 
\begin{align}\label{ope_hom}\hom(u,v)=\bigvee\{w\in\V \ | \ w\otimes u\leq v\}.\end{align}
In $\Top\cong\Cat{(\U,\two)}$, with $\xi=e_{_{\two}}^{\circ}\colon U\two\to\two$ \cite{MR2355608}, $(\V,\hom_{_{\xi}})$ is the usual {\it Sierpi\'nski space} $\mathbb{S}=(\{\bot,\top\},\{\emptyset,\{\bot\},\{\bot,\top\}\})$, and by analogy we call $(\V,\hom_{_{\xi}})$ the {\it Sierpi\'nski $(\T,\V)$-space}. For $\Cat{\V}$ ($\T=\I$), $\xi=1_{_{\V}}\colon\V\to\V$ \cite{MR2355608}, hence $\hom_{_{\xi}}=\hom$, and the Sierpi\'nski $\V$-space has easy descriptions for our examples of quantales in Table (\ref{tab1}): for $\V=\two$, it is given by the ordered set $(\{\bot,\top\},\leq)$ with $\bot<\top$; for $\V=\categ{P}_{_{+}}$, $\V=\categ{P}_{_{{\rm max}}}$, and $\V=[0,1]_{_{\odot}}$ it is given, respectively, by $([0,\infty],\ominus)$, $([0,\infty],\ovee)$, and $([0,1],\circledast)$, where, for each $u,v\in[0,\infty]$, 
\begin{align}\label{def_sierp_met}v\ominus u=\left\{\begin{array}{l}v-u, \ \text{if} \ v\geq u \\ 0, \ \text{otherwise}\end{array}\right.\qquad\with\qquad u\ovee v=\left\{\begin{array}{l}v, \ \text{if} \ u<v \\ 0, \ \text{otherwise,}\end{array}\right.\end{align}
and for each $u',v'\in[0,1]$, $u'\circledast v'={\rm min}(1,1-u+v)$. Moreover, since we are assuming $T\one=\one$, by \cite[Lemma 4.18]{MR2640215} (see also the comment after the proof of \cite[Theorem 2.9]{MR2729224}), $(\V,\hom_{_{\xi}})$ is an injective $(\T,\V)$-space and, consequently, so is the binary product $(\V\times\V,\hom_{_{\xi}}\times\hom_{_{\xi}})$, where, for each $\q\in T(\V\times\V)$, $(u,v)\in\V\times\V$, 
\begin{align*}(\hom_{_{\xi}}\times\hom_{_{\xi}})(\q,(u,v))=\hom_{_{\xi}}(T\pi_{_{1}}(\q),u)\wedge\hom_{_{\xi}}(T\pi_{_{2}}(\q),v).\end{align*}
Since $\xi\colon T\V\to\V$ is a $\T$-algebra, for $(\V,\hom_{_{\xi}})$ the order (\ref{order_space}) gives exactly the order of $\V$, which is anti-symmetric, since $\V$ is a quantale. Hence $(\V,\hom_{_{\xi}})$ and, consequently, $(\V\times\V,\hom_{_{\xi}}\times\hom_{_{\xi}})$ are separated $(\T,\V)$-spaces, whence the equivalence $\simeq$ in diagram (\ref{diag_ext}) is an equality when considering extensions with codomain $\V\times\V$; this fact will be used in Subsection \ref{sub_alex}.

For $\V$-relations $r\colon X\rel X'$ and $s\colon Y\rel Y'$, consider the $\V$-relation $r\owedge s\colon X\times Y\rel X'\times Y'$, given by, for $(x,y)\in X\times Y$, $(x',y')\in X'\times Y'$, $r\owedge s((x,y),(x',y'))=r(x,x')\wedge s(y,y')$. Assume that the diagram
\begin{align*}\xymatrix{T(\V\times\V) \ar[rr]^(0.55){T(\wedge)} \ar[d]_{\left\langle \xi\cdot T\pi_{_{1}},\xi\cdot T\pi_{_{2}}\right\rangle} \ar@{}[rrd]|{\leq} & & TV \ar[d]^{\xi} \\ \V\times\V \ar[rr]_(0.55){\wedge} & & \V}\end{align*}
is lax commutative, what is true for the examples in Table (\ref{tab1}) \cite{MR2355608, MR3227304}. Thus, by \cite[7.4]{MR3987969}, for all $\V$-relations $r\colon X\to X'$ and $s\colon Y\to Y'$, the following diagram is lax commutative,
\begin{align}\label{diag_can}\xymatrix{T(X\times Y) \ar[rr]^{{\rm can}_{_{X,Y}}} \ar[d]|{\object@{|}}_{T(r\owedge s)} \ar@{}[rrd]|{\geq} & & TX\times TY \ar[d]|{\object@{|}}^{(Tr)\owedge(Ts)} \\ T(X'\times Y') \ar[rr]_{{\rm can}_{_{X',Y'}}} & & TX'\times TY'}\end{align}
where, under the previous notation, ${\rm can}_{_{X,Y}}=\left\langle T\pi_{_{X}},T\pi_{_{Y}}\right\rangle\colon T(X\times Y)\to TX\times TY$. Then, by \cite[3.1]{MR3987969}, a $(\T,\V)$-space $(X,a)$ is exponentiable provided that, for each $\xx\in TTX$, $x\in X$, $u,v\in\V$,
\begin{align}\label{expo}\DS\bigvee_{_{\x\in TX}}(Ta(\xx,\x)\wedge u)\otimes(a(\x,x)\wedge v)\geq a(m_{_{X}}(\xx),x)\wedge(u\otimes v).\end{align}
 
Now, within the framework of strict topological theories, the tensor product of $\V$ induces a tensor product between $(\T,\V)$-spaces: for each $(X,a),(Y,b)$, $(X,a)\otimes(Y,b)=(X\times Y,c)$, where, for each $\w\in T(X\times Y)$, $(x,y)\in X\times Y$, $c(\w,(x,y))=a(T\pi_{_{X}}(\w),x)\otimes b(T\pi_{_{Y}}(\w),y)$ (see \cite[Lemma 6.1]{MR2355608}). Consider the maps
\begin{align}\label{inj1}\xymatrix{\V\otimes\V \ar[rr]^{\otimes} & & \V} \ \ \ \with \ \ \ \xymatrix{X \ar[rr]^{(-,u)} & & X\otimes\V,}\end{align}
and define, for a $\V$-relation $r\colon X\rel Y$ and an element $u\in\V$, the $\V$-relation $r\otimes u\colon X\rel Y$ given by, for each $(x,y)\in X\times Y$, $(r\otimes u)(x,y)=r(x,y)\otimes u$. Finally, consider the condition:
\begin{align}\label{inj3}\forall u,v,w\in\V, \ w\wedge(u\otimes v)=\{u'\otimes v' \ | \ u'\leq u, \ v'\leq v, \ u'\otimes v'\leq w\},\end{align}
which is equivalent to exponentiability of every injective $\V$-space in $\Cat{\V}$ \cite[Theorem 5.3]{MR2981702}. Then, by \cite[5.8]{MR3987969}, we have \vspace{0.2cm}

\noindent{\bf Theorem} {\it Suppose that: \\
$\bullet$ for all $\V$-relations $r\colon X\rel X'$ and $s\rel Y\to Y'$, diagram {\rm (\ref{diag_can})} is commutative; \\
$\bullet$ for every injective $(\T,\V)$-space $(X,a)$ and every $u\in\V$, the maps $\otimes$ and $(-,u)$ in {\rm (\ref{inj1})} are $(\T,\V)$-continuous, and $T(a\otimes u)=Ta\otimes u$; and \\
$\bullet$ {\rm (\ref{inj3})} holds. \\
Then every injective $(\T,\V)$-space is exponentiable in $\Cat{(\T,\V)}$.} \vspace{0.2cm}

\noindent{\bf Remark} As examples, in the categories of Table (\ref{tab1}) every injective $(\T,\V)$-space, and in particular the Sierpi\'nski $(\T,\V)$-space $(\V,\hom_{_{\xi}})$, is exponentiable.


\subsection{Compact and Hausdorff \texorpdfstring{$(\T,\V)$}{(T,V)}-spaces}\label{sub_comp_haus}

As discussed in \cite[Notes on Chapter V]{MR3307673}, the work of Manes \cite{MR0367901} can be considered a predecessor of the notions of compactness and Hausdorff separation in categories of $\V$-relational algebras. We follow here the definitions outlined in \cite[V-1.1.1]{MR3307673}: a $(\T,\V)$-space $(X,a)$ is {\it compact} if $1_{_{TX}}\leq a^{\circ}\cdot a$, or componentwise, if for each $\x\in TX$,
\begin{align*}k\leq\DS\bigvee_{_{x\in X}}a(\x,x)\otimes a(\x,x);\end{align*}
and it is {\it Hausdorff} if $a\cdot a^{\circ}\leq1_{_{X}}$, that is, for each $x,y\in X$, $\x\in TX$,
\begin{align*}(\bot<a(\x,x)\otimes a(\x,y) \ \Longrightarrow \ x=y) \ \ \with \ \ a(\x,x)\otimes a(\x,x)\leq k.\end{align*}
Under our assumption that $\V$ is integral, the second condition for Hausdorff separation holds trivially. Observe that, for $\Cat{(\U,\two)}\cong\Top$, we recover the fact that a topological space $(X,\tau)$ is compact and Hausdorff if, and only if, each ultrafilter has a unique convergence point. 

Under the condition that $\V$ is {\it lean}, that is, for each $u,v\in\V$,
\begin{align*}(u\vee v=\top \ \text{and} \ u\otimes v=\bot) \ \Longrightarrow \ (u=\top \ \text{or} \ v=\top),\end{align*}
by \cite[V-1.2.1]{MR3307673}, a $(\T,\V)$-space $(X,a)$ is compact and Hausdorff precisely when $(X,a)$ is a $\T$-algebra. Condition (\ref{expo}) is satisfied by any $\T$-algebra, whence compact Hausdorff $(\T,\V)$-spaces are exponentiable. Moreover, under those conditions {\bf --} $\V$ integral and lean, and a flat lax extension {\bf --} since 
\begin{align}\label{eq_comp_haus}\Cat{(\T,\V)}_{_{{\rm CompHaus}}}\cong\Set^{\T},\end{align}
by \cite[V-1.2.3]{MR3307673}, limits of compact Hausdorff $(\T,\V)$-spaces are compact and Hausdorff. Finally, by \cite[V-1.1.6(2)]{MR3307673}, assuming that the functor $T$ preserves finite coproducts, we have that finite coproducts of compact Haudorff $(\T,\V)$-spaces are compact and Hausdorff. Observe that the examples of Table (\ref{tab1}) satisfy the conditions discussed above (see \cite{MR894390}), so the latter facts, which we summarize below, hold for those categories: \\
(I) each compact Hausdorff $(\T,\V)$-space is exponentiable in $\Cat{(\T,\V)}$; \\
(II) the binary product of compact Haudorff $(\T,\V)$-spaces is compact and Hausdorff; \\
(III) finite coproducts of compact Haudorff $(\T,\V)$-spaces are compact and Hausdorff.\vspace{0.2cm}

\noindent{\bf Examples} \cite[V-1.1.4]{MR3307673} (1) For $\Cat{\V}$ with $\V$ integral and lean, by (\ref{eq_comp_haus}), $\Cat{\V}_{_{{\rm CompHaus}}}\cong\Set^{\I}$, hence a $\V$-space $(X,a)$ is compact and Hausdorff if, and only if, it is {\it discrete}, that is, $a=1_{_{X}}$. 

\noindent(2) For $\Top\cong\Cat{(\U,\two)}$, we already observed that compactness and Hausdorff separation coincide with the classical notions: each ultrafilter converges to a unique point. 

\noindent(3) For $\App\cong\Cat{(\U,\PP_{_{+}})}$ and $\NA\text{-}\App\cong\Cat{(\U,\PP_{_{{\rm max}}})}$, a space $(X,a)$ is compact exactly when it is {\it 0-compact}, that is, for each $\x\in UX$, $\DS\inf\{a(\x,x) \ | \ x\in X\}=0$; it is Hausdorff exactly when, for each $x,y\in X$, $\x\in UX$,
\begin{align*}
    (a(\x,x)<\infty \ \with \ a(\x,y)<\infty) \ \Longrightarrow \ (x=y),
\end{align*} 
or, equivalently, if the {\it pseudotopological modification} of $(X,a)$ is Hausdorff. As observed in \cite[V-1.2.2(1)]{MR3307673}, by (\ref{eq_comp_haus}), we have:
\begin{align*}\Cat{(\U,[0,1]_{_{\odot}})}_{_{{\rm CompHaus}}}\cong\NA\text{-}\App_{_{{\rm CompHaus}}}\cong\App_{_{{\rm CompHaus}}}\cong\Top_{_{{\rm CompHaus}}}\cong\Set^{\U}.\end{align*}


\section{\texorpdfstring{$\C$}{C}-generated \texorpdfstring{$(\T,\V)$}{(T,V)}-spaces}\label{sec2}


\subsection{The category \texorpdfstring{$\Cat{(\T,\V)}_{_{\C}}$}{(T,V)-CatC}}\label{sub_cgen}

This section should be compared with \cite[Section 3]{MR2080286}, whose work we follow directly. From now on, in order to keep the text lighter, we sometimes drop the prefix $(\T,\V)$ when referring to $(\T,\V)$-spaces, $(\T,\V)$-continuity, or $(\T,\V)$-structures; we also drop the prefix $|\text{-}|$ when referring to $|\text{-}|$-initial and $|\text{-}|$-final structures, which are taken with respect to the forgetful functor $|\text{-}|\colon\Cat{(\T,\V)}\to\Set$.

We fix a class $\C\subseteq\Cat{(\T,\V)}$ of objects, containing at least one non-empty element. Although the class $\C$ is arbitrary, the reader can keep it in mind as being the class of compact Hausdorff spaces. \vspace{0.2cm}
  
\noindent{\bf Definition} The elements of $\C$ are called {\it generating spaces}. For a space $(X,a)$, a continuous map from a generating space to $(X,a)$ is called a {\it probe over $(X,a)$}, or simply a {\it probe}. The {\it $\C$-generated structure} $a^{c}$ on $X$ is the final structure with respect to all probes over $(X,a)$. A space $(X,a)$ is {\it $\C$-generated} if $a=a^{c}$. The full subcategory of $\Cat{(\T,\V)}$ of $\C$-generated spaces is denoted by $\Cat{(\T,\V)}_{_{\C}}$. \vspace{0.2cm}

By definition, for a space $(X,a)$, a map $t\colon(X,a^{c})\to(Y,b)$ is continuous if, and only if, for each probe $p\colon C\to(X,a)$, the composite $t\cdot p$ is continuous. Hence, for each space $(X,a)$, the identity map $1_{_{X}}\colon(X,a^{c})\to(X,a)$ is continuous, that is, $a^{c}\leq a$. Moreover, each generating space is $\C$-generated, for if $(D,d)$ in $\C$, then $1_{_{D}}\colon D\to D$ is a probe, hence, by definition of final structure, it is a continuous map $1_{_{D}}\colon(D,d)\to(D,d^{c})$, thus $d\leq d^{c}$, and since $d^{c}\leq d$, we have $d=d^{c}$. \vspace{0.2cm}

\noindent{\bf Lemma} {\it For a $\C$-generated space $(X,a)$, a map $f\colon(X,a)\to(Y,b)$ is continuous if, and only if, $f\colon(X,a)\to(Y,b^{c})$ is so.} \vspace{0.2cm}

\noindent{\it Proof.} For sufficiency, we can factorize $f\colon(X,a)\to(Y,b)$ as
\begin{align*}\xymatrix{(X,a) \ar[rr]^{f} & & (Y,b^{c}) \ar[rr]^{1_{_{Y}}} & & (Y,b),}\end{align*}
which is continuous, since $1_{_{Y}}\colon(Y,b^{c})\to(Y,b)$ is continuous. For necessity, for each probe \linebreak $p\colon C\to(X,a)$, the composite $f\cdot p\colon C\to(Y,b)$ is continuous, hence a probe, then $f\cdot p\colon C\to(Y,b^{c})$ is continuous, and we conclude that $f\colon(X,a^{c})=(X,a)\to(Y,b^{c})$ is continuous.\qed\vspace{0.2cm}

\noindent{\bf Remark} {\it For each space $(X,a)$, $(X,a^{c})$ is $\C$-generated:} each probe $p\colon C\to(X,a)$ is a probe over $(X,a^{c})$, and, consequently, it is a probe over $(X,(a^{c})^{c})$. Hence $1_{_{X}}\colon(X,a^{c})\to(X,(a^{c})^{c})$ is a continuous map, that is, $a^{c}\leq(a^{c})^{c}$, and since $(a^{c})^{c}\leq a^{c}$, we conclude $(a^{c})^{c}=a^{c}$. \vspace{0.2cm}

The next result is also proven for the particular case of $\Top$ in \cite{Mac71}. \vspace{0.2cm}

\noindent{\bf Theorem} {\it $\Cat{(\T,\V)}_{_{\C}}$ is coreflective in $\Cat{(\T,\V)}$.} \vspace{0.2cm}

\noindent{\it Proof.} For each $(Y,b)$ in $\Cat{(\T,\V)}$, $(Y,b^{c})$ belongs to $\Cat{(\T,\V)}_{_{\C}}$. Moreover, each continuous map \linebreak $f\colon(X,a)\to(Y,b)$, with $(X,a)$ in $\Cat{(\T,\V)}_{_{\C}}$, factorizes through the identity $1_{_{Y}}\colon(Y,b^{c})\to(Y,b)$:
\begin{align*}\xymatrix@R=1.2em{(Y,b^{c}) \ar[rr]^{1_{_{Y}}} & & (Y,b) \\ & & \\ (X,a). \ar[uu]^{f} \ar[rruu]_(0.6){f} & & }\end{align*}
So the coreflector $G_{_{\C}}$ takes $(X,a)$ to $(X,a^{c})$ and $f\colon(X,a)\to(Y,b)$ to $f\colon(X,a^{c})\to(Y,b^{c})$; each coreflection is given by an identity. 
\begin{align*}\xymatrix{\Cat{(\T,\V)}_{_{\C}} \ar@{^{(}->}[rr]^{\top} & & \Cat{(\T,\V)} \ar@/_1.2pc/[ll]_{G_{_{\C}}}}\end{align*}\qed\vspace{0.2cm}

As a corollary, we obtain that $\Cat{(\T,\V)}_{_{\C}}$ is complete and cocomplete, since $\Cat{(\T,\V)}$ is so. Next we characterize the $\C$-generated spaces in terms of colimits. To do so, recall from the first section that each $|\text{-}|$-final lifting of a sink is actually the $|\text{-}|$-final lifting of a small sink. We use also the fact that each constant map is continuous. \vspace{0.2cm}

\noindent{\bf Proposition} {\it (1) $\C$-generated spaces are closed under the formation of coproducts and coequalizers, hence closed under colimits. 

\noindent(2) A space is $\C$-generated if, and only if, it is a coequalizer of a coproduct of generating spaces.} \vspace{0.2cm}

\noindent{\it Proof.} (1) This follows from the fact that the inclusion functor is a left adjoint, so it preserves colimits. 

\noindent(2) By the first assertion, a coequalizer of a coproduct of generating spaces is $\C$-generated, since generating spaces are $\C$-generated. Now let $(X,a)$ be a $\C$-generated space, hence $a=a^{c}$ is the final structure with respect to a sink of continuous maps $(p_{_{i}}\colon(X_{_{i}},a_{_{i}})\to(X,a))_{_{i\in I}}$, with $(X_{_{i}},a_{_{i}})$ in $\C$, where $I$ is a set. Take the coproduct $(\DS\dot{\bigcup_{_{i\in I}}}X_{_{i}},a_{_{I}})$ in $\Cat{(\T,\V)}$. From its universal property we get a continuous map $t\colon(\DS\dot{\bigcup_{_{i\in I}}}X_{_{i}},a_{_{I}})\to(X,a)$ such that, for each $i\in I$, $t\cdot\iota_{_{i}}=p_{_{i}}$, where $\iota_{_{i}}$ is the canonical inclusion of $X_{_{i}}$ into the coproduct.
\begin{align*}\xymatrix{(X_{_{i}},a_{_{i}}) \ar@{^{(}->}[rr]^{\iota_{_{i}}} \ar[drr]_{p_{_{i}}} & & (\DS\dot{\bigcup}X_{_{i}},a_{_{I}})  \ar[d]^{t} \\ & & (X,a)}\end{align*}
Let us prove that $t$ is a final surjection, hence a regular epimorphism in $\Cat{(\T,\V)}$. We can assume that all points of $X$ are covered by probes: for $x_{_{0}}\in X$, consider a constant map $x_{_{0}}\colon C_{_{0}}\to(X,a)$, for a non-empty element $C_{_{0}}$ of $\C$, which exists by our assumptions; adding those maps to our set of probes does not affect its finality neither its smallness, so we consider, without loss of generality, that those constant probes are already indexed by $I$. Hence, for each $x\in X$, there exists $i\in I$ such that $x=p_{_{i}}(x_{_{i}})=t(\iota_{_{i}}(x_{_{i}}))$, and $t$ is surjective. Next consider a map $s\colon(X,a)\to(Y,b)$ such that $s\cdot t$ is a continuous map; this is equivalent to, for each $i\in I$, the map $s\cdot t\cdot\iota_{_{i}}$ being continuous, hence, for each $i\in I$, $s\cdot p_{_{i}}$ is continuous, what implies that $s$ is a continuous map, since the structure $a$ is final with respect to the sink $(p_{_{i}})_{_{i\in I}}$.\qed\vspace{0.2cm}

For a complete account on regular epimorphisms in $\Cat{(\T,\V)}$ see \cite{MR2134293}. As a corollary, we can conclude that $\Cat{(\T,\V)}_{_{\C}}$ is the {\it coreflective hull} of $\C$ in $\Cat{(\T,\V)}$. Hence the matter of cartesian closedness fits the goals of \cite{MR0480687}, where the author also established the conditions used below. However, following the lines of \cite{MR2080286}, a direct approach to the question is given. 


\subsection{The category \texorpdfstring{$\C\text{-}{\rm Map}$}{C-Map}}\label{sub_cmap} 

We start by the following: \vspace{0.2cm}

\noindent{\bf Definition} A map $f\colon(X,a)\to(Y,b)$ is {\it $\C$-continuous} if the composite $f\cdot p\colon C\to(Y,b)$ is continuous, for every probe $p\colon C\to(X,a)$.\vspace{0.2cm}

Notice that, for spaces $(X,a),(Y,b)$ and a map $f\colon X\to Y$, the following assertions are equivalent: \\
(i) $f\colon(X,a)\to(Y,b)$ is $\C$-continuous; \\
(ii) $f\colon(X,a^{c})\to(Y,b)$ is continuous; \\
(iii) $f\colon(X,a^{c})\to(Y,b^{c})$ is continuous. \\
Continuity obviously implies $\C$-continuity, and from (ii) we see that for maps defined on $\C$-generated spaces the converse is also true.\vspace{0.2cm}

\noindent{\bf Lemma} {\it (1) $(\T,\V)$-spaces and $\C$-continuous maps form a category, denoted by $\C\text{-}{\rm Map}$. 

\noindent(2) The identity map $1_{_{X}}\colon(X,a^{c})\to(X,a)$ is an isomorphism in $\C\text{-}{\rm Map}$. 

\noindent(3) The assignment that sends a space $(X,a)$ to $(X,a^{c})$, and a $\C$-continuous map to itself, is an equivalence of categories $F_{_{\C}}\colon\C\text{-}{\rm Map}\to\Cat{(\T,\V)}_{_{\C}}$.} \vspace{0.2cm}

\noindent{\it Proof.} (1) Identity maps and composition of $\C$-continuous maps are readily seen to be $\C$-continuous. 

\noindent(2) The identity map $1_{_{X}}\colon(X,a^{c})\to(X,a)$ is continuous, hence it is $\C$-continuous, and \linebreak $1_{_{X}}\colon(X,a^{c})\to(X,a^{c})$ is continuous, thus $1_{_{X}}\colon(X,a)\to(X,a^{c})$ is also $\C$-continuous. 

\noindent(3) $\Cat{(\T,\V)}_{_{\C}}$ is a full subcategory of $\C\text{-}{\rm Map}$ and the inclusion $\Cat{(\T,\V)}_{_{\C}}\hookrightarrow\C\text{-}{\rm Map}$ is essentially surjective: for each space $(X,a)$, by item (2), $(X,a^{c})\cong(X,a)$ in $\C\text{-}{\rm Map}$. \qed \vspace{0.2cm} 

\subsection{Cartesian closedness of \texorpdfstring{$\C\text{-}{\rm Map}$}{C-Map} and \texorpdfstring{$\Cat{(\T,\V)}_{_{\C}}$}{(T,V)-CatC}}\label{sub_cart_clos}

Firstly we prove the following result: \vspace{0.2cm}

\noindent{\bf Lemma} {\it For $(\T,\V)$-spaces $(X,a),(Y,b),(Z,c)$, if $f\colon(X\times Y,a\times b)\to(Z,c)$ is a $\C$-continuous map, then, for each $x\in X$, the map $f_{_{x}}\colon(Y,b)\to(Z,c)$, given by $f_{_{x}}(y)=f(x,y)$ for each $y\in Y$, is $\C$-continuous.}\vspace{0.2cm}

\noindent{\it Proof.} Every constant map $x\colon(Y,b)\to(X,a)$ is continuous. Then $\left\langle x,1_{_{Y}}\right\rangle\colon(Y,b)\to(X\times Y,a\times b)$ is a continuous map, hence it is $\C$-continuous, and so is the composite $f_{_{x}}=f\cdot\left\langle x,1_{_{Y}}\right\rangle$.\qed\vspace{0.2cm}

This provides, for each $\C$-continuous map $f\colon X\times Y\to Z$, a map $\overline{f}\colon X\to\C\text{-}{\rm Map}(Y,Z)$ given by, for each $x\in X$, $\overline{f}(x)=f_{_{x}}$; as usual we call $\overline{f}$ the {\it transpose} of $f$. We wish to endow $\C\text{-}{\rm Map}(Y,Z)$ with a $(\T,\V)$-structure $d$ such that $f$ is $\C$-continuous if, and only if, $\overline{f}$ is so. In order to do that, we assume the condition: \vspace{0.2cm}

\noindent{\bf (EP)} {\it each element of $\C$ is exponentiable in $\Cat{(\T,\V)}$ and the product of two elements of $\C$ is a $\C$-generated space.} \vspace{0.2cm}

\noindent The class $\C$ is referred as being {\it productive} \cite[Definition 3.5]{MR2080286} (see also \cite{MR0480687, MR0296126}). 

Consider the spaces $(Y,b),(Z,c)$ and the sink $(q_{_{j}}\colon(Y_{_{j}},b_{_{j}})\to(Y,b))_{_{j\in J}}$ of all probes over $(Y,b)$. Since each generating space $(Y_{_{j}},b_{_{j}})$ is exponentiable, we have an exponential $(Z^{Y_{_{j}}},d_{_{j}})$ in $\Cat{(\T,\V)}$, which is given by the set
\begin{align*}Z^{Y_{_{j}}}=\{h\colon(Y_{_{j}},b_{_{j}})\to(Z,c) \ | \ h \ \text{is a $(\T,\V)$-continuous map}\}\end{align*}
endowed with the following $(\T,\V)$-structure: for each $\p\in T(Z^{Y_{_{j}}})$, $h\in Z^{Y_{_{j}}}$,
\begin{align*}d_{_{j}}(\p,h)=\bigvee\{v\in\V \ | \ (\forall\q\in (T\pi_{_{Z}})^{-1}(\p)) \ (\forall y_{_{j}}\in Y_{_{j}}) \ b_{_{j}}(T\pi_{_{X}}(\q),y_{_{j}})\wedge v\leq c(T{\rm ev}(\q),h(x))\},\end{align*}
where $\pi_{_{X}}$ and $\pi_{_{Z}}$ are the product projections from $X\times Z$ into $X$ and $Z$, respectively \cite{CHT03}. Consequently, each probe $q_{_{j}}\colon(Y_{_{j}},b_{_{j}})\to(Y,b)$ induces a map 
\begin{align*}\begin{array}{rcl} t_{q_{_{j}}}\colon\C\text{-}{\rm Map}(Y,Z) & \longrightarrow & (Z^{Y_{_{j}}},d_{_{j}}) \\ g & \longmapsto & g\cdot q_{_{j}},\end{array}\end{align*} 
which is well-defined: if $g$ is $\C$-continuous, then $g\cdot q_{_{j}}$ is continuous. We endow $\C\text{-}{\rm Map}(Y,Z)$ with the initial structure $d$ with respect to the source $(t_{q_{_{j}}}\colon\C\text{-}{\rm Map}(Y,Z)\to(Z^{Y_{_{j}}},d_{_{j}}))_{_{j\in J}}$, so that \linebreak $d=\DS\bigwedge_{_{j\in J}}t_{_{q_{_{j}}}}^{\circ}\cdot d_{_{j}}\cdot Tt_{_{q_{_{j}}}}$ is such that a map $h\colon(W,l)\to(\C\text{-}{\rm Map}(Y,Z),d)$, with $(W,l)$ in $\Cat{(\T,\V)}$, is continuous precisely when, for each $j\in J$, $t_{_{q_{_{j}}}}\cdot h\colon(W,l)\to(Z^{Y_{_{j}}},d_{_{j}})$ is continuous.\vspace{0.2cm}

\noindent{\bf Proposition} {\it For $(\T,\V)$-spaces $(X,a),(Y,b),(Z,c)$, a map $f\colon(X\times Y,a\times b)\to(Z,c)$ is $\C$-continuous if, and only if, $\overline{f}\colon(X,a)\to(\C\text{-}{\rm Map}(Y,Z),d)$ is $\C$-continuous.} \vspace{0.2cm}

\noindent{\it Proof.} Suppose that $f\colon X\times Y\to Z$ is $\C$-continuous. To prove that $\overline{f}\colon X\to\C\text{-}{\rm Map}(Y,Z)$ is $\C$-continuous, take a probe $p\colon C\to(X,a)$ and consider the composite $\overline{f}\cdot p\colon C\to\C\text{-}{\rm Map}(Y,Z)$, which we wish to verify to be a continuous map. By definition of $d$, it suffices to prove that, for each probe $q_{j}\colon(Y_{_{j}},b_{_{j}})\to(Y,b)$, $t_{_{q_{_{j}}}}\cdot\overline{f}\cdot p\colon C\to Z^{Y_{_{j}}}$ is a continuous map. We have a natural bijection
\begin{align}\label{bijection}\Cat{(\T,\V)}(C,Z^{Y_{_{j}}})\cong\Cat{(\T,\V)}(C\times Y_{_{j}},Z),\end{align}
and we calculate, for each $c\in C$, $y_{_{j}}\in Y_{_{j}}$, 
\begin{align*}(t_{_{q_{_{j}}}}\cdot\overline{f}\cdot p(c))(y_{_{j}})=\overline{f}\cdot p(c)(q_{_{j}}(y_{_{j}}))=f(p(c),q_{_{j}}(y_{_{j}}))=f\cdot(p\times q_{_{j}})(c,y_{_{j}}),\end{align*}
hence $t_{_{q_{_{j}}}}\cdot\overline{f}\cdot p$ corresponds to the map $f\cdot(p\times q_{_{j}})\colon C\times Y_{_{j}}\to Z$ by the bijection (\ref{bijection}), which is continuous, since $p\times q_{_{j}}$ is continuous, $C\times Y_{_{j}}$ is $\C$-generated, and, by hypothesis, $f$ is $\C$-continuous.

Now let $\overline{f}\colon(X,a)\to(\C\text{-}{\rm Map}(Y,Z),d)$ be $\C$-continuous. To prove that $f\colon X\times Y\to Z$ is $\C$-continuous, take a probe $r\colon C\to(X\times Y,a\times b)$ and consider the composite $f\cdot r\colon C\to Z$. Composing with the product projections $\pi_{_{X}}$ and $\pi_{_{Y}}$, we get the probes $r_{_{X}}=\pi_{_{X}}\cdot r\colon C\to X$ and \linebreak $r_{_{Y}}=\pi_{_{Y}}\cdot r\colon C\to Y$. By hypothesis, $\overline{f}\cdot r_{_{X}}\colon C\to\C\text{-}{\rm Map}(Y,Z)$ is a continuous map, what implies, by definition of $(\C\text{-}{\rm Map}(Y,Z),d)$, that $t_{r_{_{Y}}}\cdot\overline{f}\cdot r_{_{X}}\colon C\to Z^{C}$ is a continuous map. Then, for each $c\in C$, we have
\begin{align*}f\cdot r(c)=f\cdot\left\langle r_{_{X}},r_{_{Y}}\right\rangle(c)=\overline{f}(r_{_{X}}(c))(r_{_{Y}}(c))=(t_{r_{_{Y}}}\cdot\overline{f}\cdot r_{_{X}}(c))(c),\end{align*}
and we conclude that $f\cdot r$ is a continuous map.\qed\vspace{0.2cm}

\noindent{\bf Corollary} {\it $\C\text{-}{\rm Map}$ is a cartesian closed category.} \vspace{0.2cm}

\noindent{\it Proof.} For spaces $(Y,b),(Z,c)$, the evaluation map 
\begin{align*}\begin{array}{rcl}{\rm ev}_{_{Y,Z}}\colon(\C\text{-}{\rm Map}(Y,Z)\times Y,d\times b) & \longrightarrow & (Z,c) \\ (f,y) & \longmapsto & f(y),\end{array}\end{align*}
is $\C$-continuous, since its transpose $\overline{{\rm ev}}_{_{Y,Z}}\colon\C\text{-}{\rm Map}(Y,Z)\to\C\text{-}{\rm Map}(Y,Z)$ is an identity map, hence ($\C$-)continuous. Moreover, for each $\C$-continuous map $f\colon(X\times Y,a\times b)\to(Z,c)$, the $\C$-continuous map $\overline{f}\colon(X,a)\to(\C\text{-}{\rm Map}(Y,Z),d)$ is the unique one such that ${\rm ev}_{_{Y,Z}}\cdot(\overline{f}\times1_{_{Y}})=f$.
\begin{align*}\xymatrix{\C\text{-}{\rm Map}(Y,Z) & \C\text{-}{\rm Map}(Y,Z)\times Y \ar[rr]^{\hspace{1cm}{\rm ev}_{_{Y,Z}}} & & Z \\ X \ar@{.>}[u]^{\exists \ ! \ \overline{f} \ } & X\times Y \ar[u]^{\overline{f}\times 1_{_{Y}}} \ar[urr]_(0.6){f} & & }\end{align*} \qed \vspace{0.2cm}

By Lemma \ref{sub_cmap}(3) and the previous corollary, we conclude: \vspace{0.2cm}

\noindent{\bf Theorem} {\it $\Cat{(\T,\V)}_{_{\C}}$ is a cartesian closed category.} \vspace{0.2cm}

The exponential of objects $(X,a),(Y,b)$ in $\Cat{(\T,\V)}_{_{\C}}$ is given by 
\begin{align*}\C(\C\text{-}{\rm Map}((X,a),(Y,b)),d)=(\C\text{-}{\rm Map}((X,a),(Y,b)),d^{c})=(\Cat{(\T,\V)}((X,a),(Y,b)),d^{c}).\end{align*}


\section{Examples of \texorpdfstring{$\C$}{C}-generated \texorpdfstring{$(\T,\V)$}{(T,V)}-spaces}


\subsection{Compactly generated \texorpdfstring{$(\T,\V)$}{(T,V)}-spaces} 

Let $\C\subseteq\Cat{(\T,\V)}$ be the class of compact Hausdorff spaces; $\C$-generated spaces are called, as usual, {\it compactly generated}. As discussed in Subsection \ref{sub_comp_haus}, $\C$ satisfies condition {\bf (EP)}, so Theorem \ref{sub_cart_clos} holds. Together with Proposition \ref{sub_cgen}(3), we conclude that each compactly generated space is a coequalizer of a coproduct of compact Hausdorff spaces, and the category $\Cat{(\T,\V)}_{_{\C}}$ of compactly generated spaces is cartesian closed.

Let us consider the categories of Table (\ref{tab1}). For $\Cat{\V}$, quotients and coproducts of discrete objects are discrete, so that $\Cat{\V}_{_{\C}}\cong\Set$. 

In $\Top$, a space $(X,a)$ belongs to $\Cat{(\U,\two)}_{_{\C}}$ if, and only if, it is a coequalizer of a coproduct of compact Hausdorff $(\U,\two)$-spaces. Then we recover the fact that a topological space is compactly generated if, and only if, it is a quotient of a disjoint sum of compact Hausdorff spaces. This is equivalent to being a quotient of a locally compact Hausdorff space. Furthermore, the category $\Top_{_{\C}}$ of compactly generated spaces and continuous functions is cartesian closed (see \cite{MR2080286}). Classical examples of compactly generated topological spaces are sequential spaces, topological manifolds and $CW$-complexes.

For (non-Archimedean) approach spaces, we recall the equivalences given by (\ref{eq_comp_haus}): 
\begin{align*}\Cat{(\U,\categ{P}_{_{{\rm max}}})}_{_{{\rm CompHaus}}}\cong\Cat{(\U,\categ{P}_{_{{+}}})}_{_{{\rm CompHaus}}}\cong\Set^{\U}.\end{align*}
The embedding of $\Top$ in $\App$ corestricts to an embedding into $\NA$-$\App$, and $\Top$ is coreflective in both categories \cite{MR1472024, MR3731477}. Therefore, compactly generated (non-Archimedean) approach spaces, which are the colimits of 0-compact Hausdorff approach spaces, are precisely the topological approach spaces induced by a compactly generated topological space. Furthermore, they form a cartesian closed category $\App_{_{\C}}=\NA\text{-}\App_{_{\C}}$. 

Let us consider $\T=\U$ and $\V=[0,1]_{_{\odot}}$. The quantale homomorphism $\iota\colon\two\to[0,1]_{_{\odot}}$, defined by $\iota(\bot)=0$ and $\iota(\top)=1$, which is {\it compatible} with the respective lax extensions of $\U$ to $\Rell{\two}$ and $\Rell{[0,1]_{_{\odot}}}$ (see \cite[III-3.5]{MR3307673}), induces the embedding $\Top\hookrightarrow\Cat{(\U,[0,1]_{_{\odot}})}$, where a topological space $(X,a)$ is assigned to $(X,\iota\cdot a)$, with $\iota\cdot a(\x,x)=\iota(a(\x,x))\in[0,1]$, for each $\x\in UX$, $x\in X$, and morphisms are unchanged. The homomorphism $\iota$ has a right adjoint $p\colon[0,1]_{_{\odot}}\to\two$, given by $p(1)=\top$ and $p(u)=\bot$, for $u\neq 1$. The quantale homomorphism $p$ is also compatible with the lax extensions of $\U$ to $\Rell{[0,1]_{_{\odot}}}$ and $\Rell{\two}$. Hence, by \cite[III-Proposition 3.5.1]{MR3307673}, the adjunction $\iota\dashv p$ induces an adjunction
\begin{align*}\Top\xymatrix@=8ex{\ar@<-1.5mm>@{^{(}->}[r]^{\top} & \ar@<-1.5mm>@/_{2mm}/[l]}\Cat{(\U,[0,1]_{_{\odot}})}.\end{align*}
Therefore, $\Top$ is coreflective in $\Cat{(\U,[0,1]_{_{\odot}})}$, what implies that $(\U,[0,1]_{_{\odot}})$-compactly generated spaces are $(\U,[0,1]_{_{\odot}})$-spaces induced by compactly generated topological spaces. 


\subsection{Alexandroff \texorpdfstring{$(\T,\V)$}{(T,V)}-spaces}\label{sub_alex}

Let us consider $\C$ as the singleton set containing the Sierpi\'nski $(\T,\V)$-space $(\V,\hom_{_{\xi}})$. By analogy with the case in $\Top$ \cite[Example (2) of Section 3]{MR2080286}, we call the $\C$-generated spaces, or $(\V,\hom_{_{\xi}})$-generated spaces, {\it Alexandroff} spaces. Hence a $(\T,\V)$-space is Alexandroff if, and only if, it is a coequalizer of a coproduct of copies of $(\V,\hom_{_{\xi}})$.  

Next we wish to verify whether the set $\C=\{(\V,\hom_{_{\xi}})\}$ satisfies condition {\bf (EP)}. When $T\one=\one$, $(\V,\hom_{_{\xi}})$ is an injective space, and, consequently, under the hypotheses of Theorem \ref{sub_inj_exp}, it is exponentiable, and this is the case for the categories in Table (\ref{tab1}) as observed in Remark \ref{sub_inj_exp}. Let us consider the binary product $(\V\times\V,\hom_{_{\xi}}\times\hom_{_{\xi}})$; we wish to verify whether this product is $(\V,\hom_{_{\xi}})$-generated. \vspace{0.2cm} 

\noindent{\bf Lemma 1} {\it If $\T=\I$ and $\V$ is integral and totally ordered, then $(\V\times\V,\hom\times\hom)$ is an Alexandroff $\V$-space.}\vspace{0.2cm}

\noindent{\it Proof.} Let us set $d=\hom\times\hom$; denoting by $d^{c}$ the Alexandroff structure on $\V\times\V$, we have $d^{c}\leq d$. Now let $(u,v),(u',v')\in\V\times\V$; we wish to verify that $d^{c}((u,v),(u',v'))\geq d((u,v),(u',v'))$. Consider the following cases: \\
\fbox{$u\leq u'$} Since the quantale is integral, $\top\otimes u=u\leq u'$, what is equivalent to $\top\leq\hom(u,u')$, whence
\begin{align*}
    d((u,v),(u',v'))=\hom(u,u')\wedge\hom(v,v')=\top\wedge\hom(v,v')=\hom(v,v').
\end{align*}
Define the maps
\begin{align*}\begin{array}{rclcccrcl}f_{_{u}}\colon\V & \longrightarrow & \V\times\V & & \text{and} & & f_{_{v'}}\colon\V & \longrightarrow & \V\times\V \\ z & \longmapsto & (u,z) & & & & z & \longmapsto & (z,v')\end{array}\end{align*}
which are continuous, since constant maps are so. Hence
\begin{align*}d^{c}((u,v),(u,v'))=d^{c}(f_{_{u}}(v),f_{_{u}}(v'))\geq\hom(v,v'),\end{align*}
and also $d^{c}((u,v'),(u',v'))\geq\hom(u,u')=\top$. From transitivity of $d^{c}$ it follows
\begin{align*}d^{c}((u,v),(u',v'))\geq d^{c}((u,v),(u,v'))\otimes d^{c}((u,v'),(u',v'))\geq\hom(v,v')\otimes\top=d((u,v),(u',v')).\end{align*}
The case $v\leq v'$ is analogous. \\
\fbox{$u>u' \ \with \ v>v'$} Let us set
\begin{align*}
    \gamma=d((u,v),(u',v'))=\hom(u,u')\wedge\hom(v,v'),
\end{align*}
and observe that, similar to what we have in the first case, $\hom(u',u)\wedge\hom(v',v)=\top$. Consider the subset $\{\gamma,\top\}\subseteq\V$ endowed with the subspace $\V$-structure. Define the map 
\begin{align*}f\colon\{\gamma,\top\}\to\V\times\V,\quad\gamma\mapsto(u',v'),\quad\top\mapsto(u,v).\end{align*} 
Since $\gamma\leq\top$, we have $\hom(\gamma,\top)=\top$, whence
\begin{align*}
    \hom(\gamma,\top)=\hom(u',u)\wedge\hom(v',v)=d(f(\gamma),f(\top));
\end{align*}
also, by formula (\ref{ope_hom}), $\hom(\top,\gamma)=\bigvee\{w\in\V \ | \ w\otimes \top=w\leq\gamma\}=\gamma$, whence
\begin{align*}
    \hom(\top,\gamma)=\hom(u,u')\wedge\hom(v,v')=d(f(\top),f(\gamma)).
\end{align*}
Thus $f$ is fully faithful, and since $(\V\times\V,\hom\times\hom)$ is a separated injective space, there exists a continuous map $\hat{f}\colon\V\to\V\times\V$ extending $f$ along the embedding of $\{\gamma,\top\}$ into $\V$:
\begin{align*}\xymatrix{\{\gamma,\top\} \ar@{^{(}->}[rr] \ar[rd]_{f} & & \V \ar[ld]^{\hat{f}} \\ & \V\times\V. & }\end{align*}
Hence $d^{c}((u,v),(u',v'))=d^{c}(\hat{f}(\top),\hat{f}(\gamma))\geq\hom(\top,\gamma)=\gamma=d((u,v),(u',v'))$.\qed\vspace{0.2cm}

Therefore, for $\T=\I$ and $\V$ integral and totally ordered, $\C$ satisfies condition {\bf (EP)}, so that Alexandroff spaces form a cartesian closed subcategory of $\Cat{\V}$. It is straightforward to verify that in $\Ord$ every space is Alexandroff. For $\Met$, $\UltMet$, and $\categ{B}_{_{1}}\Met$, Alexandroff spaces are the coequalizers of coproducts of copies of $([0,\infty],\ominus)$, $([0,\infty],\ovee)$, and $([0,1],\circledast)$, respectively (see (\ref{def_sierp_met})).

For $\Cat{(\U,\two)}\cong\Top$, as observed in \cite[Section 3]{MR2080286}, a topological space is Alexandroff in our sense if, and only if, it is Alexandroff in the classical sense, that is, if arbitrary intersections of open sets are open, which in turn is equivalent to each point to have a smallest open set containing it. This property trivially holds for the binary product $\mathbb{S}\times\mathbb{S}$ of Sierpi\'nski spaces, since its topology is finite, whence $\mathbb{S}\times\mathbb{S}$ is Alexandroff. We recover the fact that the subcategory of Alexandroff topological spaces is cartesian closed. 

In fact, it is well-known that the subcategory of Alexandroff topological spaces is equivalent to $\Ord$ (see, for instance \cite[II-5.10.5]{MR3307673}). Motivated by this case, let us consider the pair of adjoint functors (see \cite{MR1957813}, \cite[III-3.4,3.6]{MR3307673})
\begin{align}\label{adj_alex}\xymatrix{\Cat{\V} \ar@<1.0ex>[rr]^(0.45){A^{\circ}} \ar@{}[rr]|(0.45){\bot} \ar[rd]_{{\rm top.}} & & \Cat{(\T,\V)} \ar@<1.0ex>[ll]^(0.55){A_{_{e}}} \ar[ld]^{{\rm top.}} \\ & \Set, & }\end{align}
where, for each $(\T,\V)$-space $(X,a)$, $A_{_{e}}(X,a)=(X,a\cdot e_{_{X}})$, with $e_{_{X}}\colon X\to TX$ the $X$-component of the natural transformation $e\colon{\rm Id}_{_{\Set}}\to T$, and, for each $\V$-space $(Y,\beta)$, $A^{\circ}(Y,\beta)=(Y,\beta_{_{\#}})$, with $\beta_{_{\#}}=e_{_{Y}}^{\circ}\cdot T\beta$; on morphisms both functors are the identity, and the functors into $\Set$ are the forgetful functors which are topological.

When $\T=\U$ and $\V=\two$, the instance of adjunction (\ref{adj_alex}) gives
\begin{align*}\Ord\xymatrix@=8ex{\ar@{}[r]|{\bot}\ar@<1mm>@/^{2mm}/[r]^{A^{\circ}} & \ar@<1mm>@/^{2mm}/[l]^{A_{_{e}}}}\Top,\end{align*}
where to each ordered set $(X,\leq)$ is assigned the space $(X,\tau_{_{\leq}})$, with $\tau_{_{\leq}}$ the Alexandroff topology, that is, the topology that has as a basis the sets $\downarrow\! x$, $x\in X$; and to each topological space $(X,\tau)$ is assigned the ordered set $(X,\leq_{_{\tau}})$, where $\leq_{_{\tau}}$ is the dual of the {\it specialization order}, that is, 
\begin{align*}x\leq_{_{\tau}}y \ \Longleftrightarrow \ \dot{x}\rightarrow y,\end{align*} 
where $\dot{x}$ is the principal ultrafilter generated by $\{x\}$ and $\rightarrow$ denotes the convergence relation between ultrafilters and points defined by $\tau$. 

Moreover, Alexandroff topological spaces are precisely the spaces that are the image by $A^{\circ}$ of an ordered set \cite[II-5.10.5, III-3.4.3(1)]{MR3307673}. We wish to find conditions under which this fact holds in our general setting. Firstly, we must have that $(\V,\hom_{_{\xi}})$ itself is the image by $A^{\circ}$ of some Alexandroff $\V$-space. The natural candidate is $(\V,\hom)$, so we wish to verify under which conditions we have $\hom_{_{\xi}}=\hom_{_{\#}}=e_{_{\V}}^{\circ}\cdot T\hom$. By (\ref{extension}), for each $\v\in T\V$, $v\in\V$,
\begin{align*}\begin{array}{rcl}e_{_{\V}}^{\circ}\cdot T\hom(\v,v) & = & T\hom(\v,e_{_{\V}}(v)) \\[0.1cm]
& = & \bigvee\{\xi\cdot T\overrightarrow{\hom}(\w) \ | \ \w\in T(\V\times\V), \ T\pi_{_{1}}(\w)=\v, \ T\pi_{_{2}}(\w)=e_{_{\V}}(v)\}.\end{array}\end{align*}
Furthermore, by \cite[Lemma 3.2]{MR2355608}, lax commutativity of the diagram
\begin{align*}\xymatrix{T(\V\times\V) \ar[rr]^{T\overrightarrow{\hom}} \ar[d]_{\left\langle\xi\cdot T\pi_{_{1}},\xi\cdot T\pi_{_{2}}\right\rangle} \ar@{}[rrd]|{\geq} & & T\V \ar[d]^{\xi} \\ \V\times\V \ar[rr]_{\overrightarrow{\hom}} & & \V}\end{align*}
is assured. Hence, for each $\w\in T(\V\times\V)$ such that $T\pi_{_{1}}(\w)=\v$ and $T\pi_{_{2}}(\w)=e_{_{\V}}(v)$, we have $\xi\cdot T\overrightarrow{\hom}(\w)\leq\hom_{_{\xi}}(\v,v)$, and consequently $e_{_{\V}}^{\circ}\cdot T\hom(\v,v)\leq\hom_{_{\xi}}(\v,v)$. Then we can see that the required condition is strict commutativity of the latter diagram.\vspace{0.2cm} 

\noindent{\bf Theorem} {\it If the diagram below is commutative, then the functor $A^{\circ}$ preserves Alexandroff spaces.
\begin{align}\label{diag_hom_alex}\xymatrix{T(\V\times\V) \ar[rr]^{T\overrightarrow{\hom}} \ar[d]_{\left\langle\xi\cdot T\pi_{_{1}},\xi\cdot T\pi_{_{2}}\right\rangle} & & T\V \ar[d]^{\xi} \\ \V\times\V \ar[rr]_{\overrightarrow{\hom}} & & \V}\end{align}}\vspace{0.2cm}

\noindent{\it Proof.} Commutativity of (\ref{diag_hom_alex}) implies that $A^{\circ}(\V,\hom)=(\V,\hom_{_{\xi}})$. Let $(X,\alpha)$ be an Alexandroff $\V$-space, and $(X,a)=A^{\circ}(X,\alpha)$. Let $h\colon (X,a)\to(Y,b)$ be a map such that, for every continuous map $f\colon(\V,\hom_{_{\xi}})\to(X,a)$, the composite $h\cdot f$ is continuous. We wish to prove that $h$ is a $(\T,\V)$-continuous map. Since $A^{\circ}\dashv A_{_{e}}$, we only need to verify that $h\colon(X,\alpha)\to A_{_{e}}(Y,b)=(Y,b\cdot e_{_{Y}})$ is a $\V$-continuous map, which holds if, and only if, for each $\V$-continuous map $f\colon(\V,\hom)\to(X,\alpha)$, the composite $h\cdot f\colon(\V,\hom)\to(Y,b\cdot e_{_{Y}})$ is $\V$-continuous, since $(X,\alpha)$ is Alexandroff. 

Each $\V$-continuous map $f$ from $(\V,\hom)$ to $(X,\alpha)$ becomes a $(\T,\V)$-continuous map from $(\V,\hom_{_{\xi}})$ to $(X,a)$ by applying the functor $A^{\circ}$. Therefore, $h\cdot f\colon(\V,\hom_{_{\xi}})\to(Y,b)$ is a $(\T,\V)$-continuous map, whence $h\cdot f\colon(\V,\hom)\to(Y,b\cdot e_{_{Y}})$ is a $\V$-continuous map.\qed\vspace{0.2cm}  

\noindent{\bf Proposition} {\it (1) If $(X,a)$ is an Alexandroff $(\T,\V)$-space, then $(X,a)=A^{\circ}\cdot A_{_{e}}(X,a)$. \\
(2) If $\T$ is such that, for each set $X$, 
\begin{align}\label{ineq_prod_alex}(e_{_{X}}\times e_{_{X}})^{\circ}\cdot{\rm can}_{_{X,X}}\leq e_{_{X\times X}}^{\circ},\end{align}
where ${\rm can}_{_{X,X}}=\langle T\pi_{_{1}}^{X},T\pi_{_{2}}^{X}\rangle\colon T(X\times X)\to TX\times TX$, then, for each $\V$-space $(X,\alpha)$, we have $(X,\alpha)=A_{_{e}}\cdot A^{\circ}(X,\alpha)$.}\vspace{0.2cm} 

\noindent{\it Proof.} (1) Let us verify that $a=e_{_{X}}^{\circ}\cdot Ta\cdot Te_{_{X}}$. The equality $m_{_{X}}\cdot Te_{_{X}}=1_{_{TX}}$ implies the inequality $Te_{_{X}}\leq m_{_{X}}^{\circ}$, and $1_{_{X}}\leq a\cdot e_{_{X}}$ is equivalent to $e_{_{X}}^{\circ}\leq a$, whence
\begin{align*}e_{_{X}}^{\circ}\cdot Ta\cdot Te_{_{X}}\leq a\cdot Ta\cdot m_{_{X}}^{\circ}\leq a\cdot m_{_{X}}\cdot m_{_{X}}^{\circ}\leq a\end{align*}
(see \cite[III-3.4.2]{MR3307673}). For the converse inequality, by the adjunction $A^{\circ}\dashv A_{_{e}}$, each continuous map $f\colon(\V,\hom_{_{\xi}})=A^{\circ}(\V,\hom)\to(X,a)$ is continuous from $(\V,\hom)$ to $A_{_{e}}(X,a)=(X,a\cdot e_{_{X}})$, and applying $A^{\circ}$ we obtain a continuous map $f\colon(\V,\hom_{_{\xi}})\to A^{\circ}(X,a\cdot e_{_{X}})=(X,e_{_{X}}^{\circ}\cdot Ta\cdot Te_{_{X}})$. Since $(X,a)$ is Alexandroff, the identity map $1_{_{X}}$ is continuous,
\begin{align*}\xymatrix{(\V,\hom_{_{\xi}}) \ar[rr]^{f} \ar[rrd]_{f} & & (X,a) \ar[d]^{1_{_{X}}} \\ & & (X,e_{_{X}}^{\circ}\cdot Ta\cdot Te_{_{X}})}\end{align*}
hence $a\leq e_{_{X}}^{\circ}\cdot Ta\cdot Te_{_{X}}$. \\
(2) For each set $X$, $\w\in T(X\times X)$, $(x,x')\in X\times X$,
\begin{align*}\begin{array}{rrl}e_{_{X\times X}}^{\circ}(\w,(x,x'))=k & \Longleftrightarrow & e_{_{X\times X}}(x,x')=\w \\
& \Longrightarrow & \left(T\pi_{_{1}}^{X}(\w)=e_{_{X}}(x)\quad\with\quad T\pi_{_{2}}^{X}(\w)=e_{_{X}}(x')\right) \\
& \Longleftrightarrow & {\rm can}_{_{X,X}}(\w)=e_{_{X}}\times e_{_{X}}(x,x') \\
& \Longleftrightarrow & (e_{_{X}}\times e_{_{X}})^{\circ}\cdot{\rm can}_{_{X,X}}(\w,(x,x'))={\rm can}_{_{X,X}}(\w,e_{_{X}}\times e_{_{X}}(x,x'))=k,\end{array}\end{align*}
hence $e_{_{X\times X}}^{\circ}\leq(e_{_{X}}\times e_{_{X}})^{\circ}\cdot{\rm can}_{_{X,X}}$, and if (\ref{ineq_prod_alex}) holds, then it is an equality.

For each $(x,x')\in X\times X$, if $\w\in T(X\times X)$ is such that $T\pi_{_{1}}^{X}(\w)=e_{_{X}}(x)$ and $T\pi_{_{2}}^{X}(\w)=e_{_{X}}(x')$, then $e_{_{X\times X}}(x,x')=\w$, and, for each $\V$-space $(X,\alpha)$, we calculate:
\begin{align*}\begin{array}{rcl}T\alpha(e_{_{X}}(x),e_{_{X}}(x')) & = & \bigvee\{\xi\cdot T\vec{\alpha}(\w) \ | \ \w\in T(X\times X), \ T\pi_{_{1}}^{X}(\w)=e_{_{X}}(x), \ T\pi_{_{2}}^{X}(\w)=e_{_{X}}(x')\} \\
& = & \xi\cdot T\vec{\alpha}\cdot e_{_{X\times X}}(x,x') \\
& = & \xi\cdot e_{_{\V}}\cdot\vec{\alpha}(x,x')\qquad\qquad\text{($e$ is a natural transformation)} \\
& = & \vec{\alpha}(x,x')\qquad\qquad\qquad\hspace{0.3cm}\text{(because $\xi\colon T\V\to\V$ is a $\T$-algebra).}\end{array}\end{align*}
Therefore $(X,\alpha)=(X,e_{_{X}}^{\circ}\cdot T\alpha\cdot e_{_{X}})=A_{_{e}}\cdot A^{\circ}(X,\alpha)$.\qed\vspace{0.2cm}

\noindent{\bf Corollary} {\it If the diagram {\rm (\ref{diag_hom_alex})} is commutative and $\T$ satisfies {\rm (\ref{ineq_prod_alex})}, then the Alexandroff $(\T,\V)$-spaces are precisely the images by $A^{\circ}$ of Alexandroff $\V$-spaces.} \vspace{0.2cm}

\noindent{\it Proof.} By the previous Theorem, commutativity of (\ref{diag_hom_alex}) implies that the image by $A^{\circ}$ of any Alexandroff $\V$-space is an Alexandroff $(\T,\V)$-space. 

Conversely, let $(X,a)$ be an Alexandroff $(\T,\V)$-space. Then, by item (1) of the previous Proposition, $(X,a)=A^{\circ}\cdot A_{_{e}}(X,a)$. Let us verify that $(X,\alpha)=A_{_{e}}(X,a)$ is an Alexandroff $\V$-space. Let $(Y,\beta)$ be a $\V$-space and $f\colon(X,\alpha)\to(Y,\beta)$ be a map such that, for every $\V$-continuous map $p\colon(\V,\hom)\to(X,\alpha)$, the composite $f\cdot p\colon(\V,\hom)\to(Y,\beta)$ is $\V$-continuous.
\begin{align*}
    \xymatrix{(\V,\hom) \ar[rr]^{p} \ar[rrd]_{f\cdot p} & & (X,\alpha) \ar[d]^{f} \\ & & (Y,\beta)}
\end{align*}
We wish to prove that $f$ is $\V$-continuous. Since $\T$ satisfies (\ref{ineq_prod_alex}), by item (2) of the previous Proposition, $(Y,\beta)=A_{_{e}}\cdot A^{\circ}(Y,\beta)$. Hence, by the adjunction $A^{\circ}\dashv A_{_{e}}$, $f\colon(X,\alpha)\to A_{_{e}}\cdot A^{\circ}(Y,\beta)$ is $\V$-continuous if, and only if, $f\colon A^{\circ}(X,\alpha)=(X,a)\to A^{\circ}(Y,\beta)$ is $(\T,\V)$-continuous. Consider a $(\T,\V)$-continuous map $p\colon(\V,\hom_{_{\xi}})=A^{\circ}(\V,\hom)\to(X,a)$. Then $p\colon(\V,\hom)\to(X,\alpha)$ is a $\V$-continuous map. Hence, by hypothesis, the composite $f\cdot p\colon(\V,\hom)\to(Y,\beta)=A_{_{e}}\cdot A^{\circ}(Y,\beta)$ is $\V$-continuous, and so $f\cdot p\colon A^{\circ}(\V,\hom)=(\V,\hom_{_{\xi}})\to A^{\circ}(Y,\beta)$ is $(\T,\V)$-continuous. Since $(X,a)$ is an Alexandroff $(\T,\V)$-space, we conclude that $f\colon(X,a)\to A^{\circ}(Y,\beta)$ is $(\T,\V)$-continuous. Therefore, $(X,\alpha)$ is $(\V,\hom)$-generated. \qed \vspace{0.2cm}   

\noindent{\bf Lemma 2} {\it For $\V$ integral and totally ordered, if the diagrams
\begin{align*}\xymatrix{T(\V\times\V) \ar[rr]^(0.55){T(\wedge)} \ar[d]_{\left\langle \xi\cdot T\pi_{_{1}},\xi\cdot T\pi_{_{2}}\right\rangle} \ar@{}[rrd]|{\leq} & & T\V \ar[d]^{\xi} \\ \V\times\V \ar[rr]_(0.55){\wedge} & & \V}\qquad\qquad\xymatrix{T(\V\times\V) \ar[rr]^{T\overrightarrow{\hom}} \ar[d]_{\left\langle\xi\cdot T\pi_{_{1}},\xi\cdot T\pi_{_{2}}\right\rangle} & & T\V \ar[d]^{\xi} \\ \V\times\V \ar[rr]_{\overrightarrow{\hom}} & & \V}\end{align*}
are (lax) commutative and the inequality {\rm (\ref{ineq_prod_alex})} holds for $X=\V$, i.e., $(e_{_{\V}}\times e_{_{\V}})^{\circ}\cdot{\rm can}_{_{\V,\V}}\leq e_{_{\V\times\V}}^{\circ}$, then $(\V\times\V,\hom_{_{\xi}}\times\hom_{_{\xi}})$ is an Alexandroff $(\T,\V)$-space.}\vspace{0.2cm}

\noindent{\it Proof.} By Lemma 1 and the previous Theorem, it suffices to show that 
\begin{align*}(\V\times\V,\hom_{_{\xi}}\times\hom_{_{\xi}})=A^{\circ}(\V\times\V,\hom\times\hom).\end{align*}
For each $(u,v),(z,w)\in\V\times\V$,
\begin{align*}\begin{array}{rcl}(\hom_{_{\xi}}\times\hom_{_{\xi}})\cdot e_{_{\V\times\V}}((u,v),(z,w)) & = & \hom_{_{\xi}}\times\hom_{_{\xi}}(e_{_{\V\times\V}}(u,v),(z,w)) \\[0.1cm]
& = & \hom_{_{\xi}}(e_{_{\V}}(u),z)\wedge\hom_{_{\xi}}(e_{_{\V}}(v),w) \\[0.1cm]
& = & \hom\times\hom((u,v),(z,w)).\end{array}\end{align*}
Hence $(\V\times\V,(\hom\times\hom)_{_{\#}})=A^{\circ}(\V\times\V,\hom\times\hom)=A^{\circ}\cdot A_{_{e}}(\V\times\V,\hom_{_{\xi}}\times\hom_{_{\xi}})$, so that $(\hom\times\hom)_{_{\#}}\leq\hom_{_{\xi}}\times\hom_{_{\xi}}$, since the counit of $A^{\circ}\dashv A_{_{e}}$ is an identity map. Conversely, for each $\w\in T(\V\times\V)$, $(u,v)\in\V\times\V$,
\begin{align*}\begin{array}{rcl}\hom_{_{\xi}}\times\hom_{_{\xi}}(\w,(u,v)) & = & \hom_{_{\xi}}(T\pi_{_{1}}(\w),u)\wedge\hom_{_{\xi}}(T\pi_{_{2}}(\w),v) \\[0.1cm]
& = & e_{_{\V}}^{\circ}\cdot T\hom(T\pi_{_{1}}(\w),u)\wedge e_{_{\V}}^{\circ}\cdot T\hom(T\pi_{_{2}}(\w),v) \\[0.1cm]
& = & T\hom\times T\hom(\w,e_{_{V}}\times e_{_{\V}}(u,v)) \\[0.1cm]
& = & (e_{_{\V}}\times e_{_{\V}})^{\circ}\cdot(T\hom\times T\hom)(\w,(u,v)) \\[0.1cm]
& = & (e_{_{\V}}\times e_{_{\V}})^{\circ}\cdot(T\hom\owedge T\hom)\cdot{\rm can}_{_{\V,\V}}(\w,(u,v)) \\[0.1cm]
& \leq & (e_{_{\V}}\times e_{_{\V}})^{\circ}\cdot{\rm can}_{_{\V,\V}}\cdot T(\hom\owedge\hom)(\w,(u,v))\qquad\text{(by (\ref{diag_can}))} \\[0.1cm]
& = & (e_{_{\V}}\times e_{_{\V}})^{\circ}\cdot{\rm can}_{_{\V,\V}}\cdot T(\hom\times\hom)(\w,(u,v)) \\[0.1cm]
& \leq & e_{_{\V\times\V}}^{\circ}\cdot T(\hom\times\hom)(\w,(u,v))\qquad\text{(by hypothesis)} \\[0.1cm]
& = & (\hom\times\hom)_{_{\#}}(\w,(u,v)).\end{array}\end{align*}\qed\vspace{0.2cm}

Therefore, under the conditions of Lemma 2, $\C=\{(\V,\hom_{_{\xi}})\}$ satisfies condition {\bf (EP)} and Alexandroff $(\T,\V)$-spaces form a cartesian closed subcategory of $\Cat{(\T,\V)}$. \vspace{0.2cm} 

\noindent{\bf Examples} (1) Let us verify that, for the category $\App\cong\Cat{(\U,\PP_{_{+}})}$, the conditions of Lemma 2 are satisfied. In this case, $\xi\colon U[0,\infty]\to[0,\infty]$ is given by $\xi(\v)=\inf\{u\in[0,\infty] \ | \ [0,u]\in\v\}$, for each $\v\in U[0,\infty]$ (see, for instance, \cite{MR3227304}). Moreover, with the same definition of $\xi$, $(\U,\PP_{_{{\rm max}}},\xi)$ is a strict topological theory, hence the commutativity of the diagram 
\begin{align*}\xymatrix{U([0,\infty]\times[0,\infty]) \ar[rr]^(0.6){U({\rm max})} \ar[d]_{\left\langle\xi\cdot U\pi_{_{1}},\xi\cdot U\pi_{_{2}}\right\rangle} & & U[0,\infty] \ar[d]^{\xi} \\ [0,\infty]\times[0,\infty] \ar[rr]_(0.6){{\rm max}} & & [0,\infty]}\end{align*}
follows from (\ref{diag_top_the}). Consider the diagram
\begin{align*}\xymatrix{U([0,\infty]\times[0,\infty]) \ar[rr]^(0.55){U\vec{\ominus}} \ar[d]_{\left\langle\xi\cdot U\pi_{_{1}},\xi\cdot U\pi_{_{2}}\right\rangle} & & U[0,\infty] \ar[d]^{\xi} \\ [0,\infty]\times[0,\infty] \ar[rr]_(0.55){\vec{\ominus}} & & [0,\infty].}\end{align*}
We know that $\vec{\ominus}\cdot\left\langle\xi\cdot U\pi_{_{1}},\xi\cdot U\pi_{_{2}}\right\rangle\leq\xi\cdot U\vec{\ominus}$ (recall that the order on $[0,\infty]$ is $\geq$). Let us assume that $\vec{\ominus}\cdot\left\langle\xi\cdot U\pi_{_{1}},\xi\cdot U\pi_{_{2}}\right\rangle<\xi\cdot U\vec{\ominus}$. Hence there exists $\w\in U([0,\infty]\times[0,\infty])$ such that, fixing $\v_{_{1}}=U\pi_{_{1}}(\w)$ and $\v_{_{2}}=U\pi_{_{2}}(\w)$, we have
\begin{align*}\xi(\v_{_{2}})\ominus\xi(\v_{_{1}})<\xi(U\vec{\ominus}(\w))={\rm inf}\{u\in[0,\infty] \ | \ [0,u]\in U\vec{\ominus}(\w)\}.\end{align*}
Here $[0,u]\in U\vec{\ominus}(\w)$ if, and only if, $(\vec{\ominus})^{-1}([0,u])\in\w$, so that $\xi(U\vec{\ominus}(\w))={\rm inf}\{u\in[0,\infty] \ | \ S_{_{u}}\in\w\}$, where the set
\begin{align*}S_{_{u}}=(\vec{\ominus})^{-1}([0,u])=\{(p,q)\in[0,\infty]\times[0,\infty] \ | \ q\ominus p\leq u\}\end{align*}
can be depicted as the gray area in the graphic below. \vspace{0.2cm}
\begin{center}
    \includegraphics[scale=0.5]{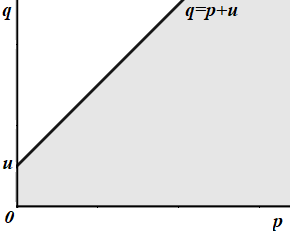}
\end{center}
Let $t\in[0,\infty]$ be such that $\xi(\v_{_{2}})\ominus\xi(\v_{_{1}})<t<\xi(U\vec{\ominus}(\w))$; then $S_{_{t}}\notin\w$. Since $\xi(\v_{_{2}})<\xi(\v_{_{1}})+t$, there exists $n\in\mathbb{N}$ such that $\xi(\v_{_{2}})+\frac{t}{n}<\xi(\v_{_{1}})+t$. Let us assume that $\xi(\v_{_{1}})>0$ so that we can choose $\frac{t}{n}<\xi(\v_{_{1}})$. Hence
\begin{align*}\xi(\v_{_{1}})-\frac{t}{n}<\xi(\v_{_{1}}) \ \Rightarrow \ [0,\xi(\v_{_{1}})-\frac{t}{n}]\notin\v_{_{1}} \ \Rightarrow \ ]\xi(\v_{_{1}})-\frac{t}{n},\infty]\in\v_{_{1}} \ \Leftrightarrow \ ]\xi(\v_{_{1}})-\frac{t}{n},\infty]\times[0,\infty]\in\w\end{align*}
and
\begin{align*}\xi(\v_{_{2}})+\frac{t}{n}>\xi(\v_{_{2}}) \ \Rightarrow \ [0,\xi(\v_{_{2}})+\frac{t}{n}]\in\v_{_{2}} \ \Leftrightarrow \ [0,\infty]\times[0,\xi(\v_{_{2}})+\frac{t}{n}]\in\w.\end{align*}
Hence $]\xi(\v_{_{1}})-\frac{t}{n},\infty]\times[0,\xi(\v_{_{2}})+\frac{t}{n}]\in\w$, but $]\xi(\v_{_{1}})-\frac{t}{n},\infty]\times[0,\xi(\v_{_{2}})+\frac{t}{n}]\subseteq S_{_{t}}$, what implies that $S_{_{t}}\in\w$, a contradiction.
\begin{center}
    \includegraphics[scale=0.5]{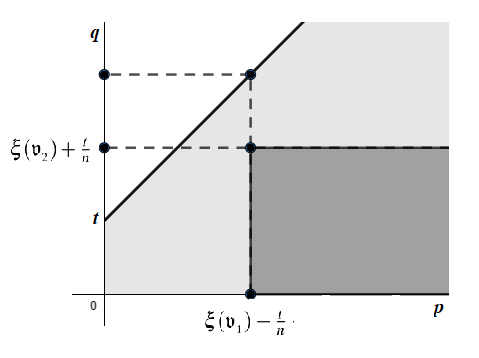}
\end{center}
In the case $\xi(\v_{_{1}})=0$, we have that $\xi(\v_{2})+\frac{t}{n}<t$, whence $[0,\infty]\times[0,\xi(\v_{2})+\frac{t}{n}]\in\w$, and $[0,\infty]\times[0,\xi(\v_{2})+\frac{t}{n}]\subseteq S_{_{t}}$, so we obtain a contradiction.

Finally, for each set $X$, if $\w\in U(X\times X)$ is such that, for $(x,x')\in X\times X$, $U\pi_{_{1}}^{X}(\w)=e_{_{X}}(x)$ and $U\pi_{_{2}}^{X}(\w)=e_{_{X}}(x')$, then $\{x\}\in U\pi_{_{1}}^{X}(\w)$ is equivalent to $(\pi_{_{1}}^{X})^{-1}(\{x\})=\{x\}\times X\in\w$, and $\{x'\}\in U\pi_{_{2}}^{X}(\w)$ is equivalent to $(\pi_{_{2}}^{X})^{-1}(\{x'\})=X\times\{x'\}\in\w$, whence 
\begin{align*}(\{x\}\times X)\cap(X\times\{x'\})=\{(x,x')\}\in\w,\end{align*} 
that is, $\w=e_{_{X\times X}}(x,x')$. Thus $(e_{_{X}}\times e_{_{X}})^{\circ}\cdot{\rm can}_{_{X,X}}\leq e_{_{X\times X}}$.

Therefore we conclude that Alexandroff approach spaces form a cartesian closed subcategory of $\App$. Moreover, Alexandroff approach spaces are precisely the images of Alexandroff metric spaces by the functor $A^{\circ}\colon\Met\to\App$, where, for each metric space $(X,d)$, $A^{\circ}(X,d)=(X,d_{_{\#}})$ with, for each $\x\in UX$, $x\in X$,
\begin{align*}d_{_{\#}}(\x,x)=\DS\sup_{A\in\x}\DS\inf_{x'\in A} d(x',x).\end{align*}
In terms of approach distances, for each $x\in X$, $A\subseteq X$,
\begin{align*}d_{_{\#}}(x,A)=\DS\inf_{x'\in A} d(x',x)\end{align*}
(see \cite[III-3.4.3(2),2.4.1(1)]{MR3307673}). Expressly, if $(X,d)$ is an Alexandroff metric space, then $(X,d_{_{\#}})$ is an Alexandroff approach space; if $(X,a)$ is an Alexandroff approach space, then $(X,a\cdot e_{_{X}})$ is an Alexandroff metric space, and, furthermore, $a=d_{_{\#}}$, with $d=a\cdot e_{_{X}}$, so that, for each $\x\in UX$, $x\in X$,
\begin{align*}a(\x,x)=\DS\sup_{A\in\x}\inf_{x'\in A} a(e_{_{X}}(x'),x).\end{align*}
In terms of approach distances, if $(X,\delta)$ is an Alexandroff approach space, then, for each $x\in X$, $A\subseteq X$, 
\begin{align*}\begin{array}{rcl}\delta(x,A) & = & \DS\inf_{x'\in A}\sup_{\scaleto{B\in e_{_{X}}(x')}{7pt}}\delta(x,B) \\[0.1cm]
& = & \DS\inf_{x'\in A}\delta(x,\{x'\}),\end{array}\end{align*}
since $\{x'\}\subseteq B$ implies $\delta(x,B)\leq\delta(x,\{x'\})$. Therefore, Alexandroff approach spaces are the metric approach spaces \cite[Theorem 3.1.11]{MR1472024}. \vspace{0.2cm}

\noindent (2) Consider the category $\Cat{(\U,[0,1]_{_{\odot}})}$. Let us verify that the conditions of Lemma 2 hold. Firstly, $\xi\colon U[0,1]\to[0,1]$ is defined by $\xi(\v)=\sup\{u\in[0,1] \ | \ [u,1]\in\v\}$. Let us verify the commutativity of the diagram
\begin{align*}\xymatrix{U([0,1]\times[0,1]) \ar[rr]^(0.55){U(\wedge)} \ar[d]_{\left\langle \xi\cdot U\pi_{_{1}},\xi\cdot U\pi_{_{2}}\right\rangle} \ar@{}[rrd]|{\leq} & & U[0,1] \ar[d]^{\xi} \\ [0,1]\times[0,1] \ar[rr]_(0.55){\wedge} & & [0,1].}\end{align*}
Let $\w\in U([0,1]\times[0,1])$ and fix $\v_{_{i}}=U\pi_{_{i}}(\w)$, $i=1,2$. Suppose that $\xi\cdot U(\wedge)(\w)<\xi(\v_{_{1}})\wedge\xi(\v_{_{2}})$. Hence there exists $t\in[0,1]$ with $\xi\cdot U(\wedge)(\w)<t<\xi(\v_{_{1}})\wedge\xi(\v_{_{2}})$. Proceeding in a similar way as in the first example, this means that $[t,1]\times[t,1]\notin\w$ and, by definition of $\xi$, $[t,1]\in\v_{_{i}}$, $i=1,2$. Then $[t,1]\times[t,1]\in\w$, a contradiction. Therefore $\xi(\v_{_{1}})\wedge\xi(\v_{_{2}})\leq\xi\cdot U(\wedge)(\w)$.   

Consider the diagram
\begin{align*}\xymatrix{U([0,1]\times[0,1]) \ar[rr]^(0.55){U\vec{\circledast}} \ar[d]_{\left\langle\xi\cdot U\pi_{_{1}},\xi\cdot U\pi_{_{2}}\right\rangle} & & U[0,1] \ar[d]^{\xi} \\ [0,1]\times[0,1] \ar[rr]_(0.55){\vec{\circledast}} & & [0,1].}\end{align*}
We know that $\xi\cdot U\vec{\circledast}\leq\vec{\circledast}\cdot\left\langle\xi\cdot U\pi_{_{1}},\xi\cdot U\pi_{_{2}}\right\rangle$. Let $\w\in U([0,1]\times[0,1])$ and suppose we have $\xi\cdot U\vec{\circledast}(\w)<\xi(\v_{_{1}})\circledast\xi(\v_{_{2}})$, where $\v_{_{i}}=U\pi_{_{i}}(\w)$, $i=1,2$. Hence there exists $t\in[0,1]$ such that $\xi\cdot U\vec{\circledast}(\w)<t<\xi(\v_{_{1}})\circledast\xi(\v_{_{2}})$. Analogously to the first item, one calculates 
\begin{align*}\xi\cdot U\vec{\circledast}(\w)=\sup\{u\in[0,1] \ | \ R_{_{u}}\in\w\},\end{align*} 
where $R_{_{u}}=\{(p,q)\in[0,1]\times[0,1] \ | \ u-1+p\leq q\}$, which is depicted as the gray area of the graphic below.
\vspace{0.2cm}
\begin{center}
    \includegraphics[scale=0.5]{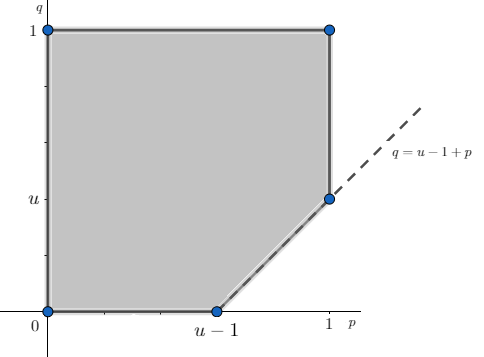}
\end{center}
Hence $R_{_{t}}\notin\w$. Let us study the cases below.

\noindent\fbox{$\xi(\v_{_{1}})=0$} Since $0<t<1$, there exists $\epsilon>0$ such that $\epsilon<1-t$. Hence
\begin{align*}
    \epsilon=\xi(\v_{_{1}})+\epsilon>\xi(\v_{_{1}}) \ \Rightarrow \ [\epsilon,1]\notin\v_{_{1}} \ \Rightarrow \ [0,\epsilon[\in\v_{_{1}} \ \Rightarrow \ [0,\epsilon[\times[0,1]\in\w,
\end{align*}
but $[0,\epsilon[\times[0,1]\subseteq R_{_{t}}$, what implies $R_{_{t}}\in\w$, a contradiction.

\noindent\fbox{$\xi(\v_{_{1}})>0$} If $\xi(\v_{_{1}})<1-t$, then we can choose $\epsilon>0$ with $\xi(\v_{_{1}})<\epsilon<1-t$. Hence $[0,\epsilon[\times[0,1]\in\w$ and $[0,\epsilon[\times[0,1]\subseteq R_{_{t}}$, a contradiction. Let us assume that $\xi(\v_{_{1}})\geq 1-t$. Since $t<\xi(\v_{_{1}})\circledast\xi(\v_{_{2}})$, then $t<1-\xi(\v_{_{1}})+\xi(\v_{_{2}})$. Choose $\epsilon>0$ such that $\xi(\v_{_{2}})-\epsilon>t-1+\xi(\v_{_{1}})\geq 0$. Thus
\begin{align*}
    \xi(\v_{_{1}})<\xi(\v_{_{2}})-\epsilon-t+1 \ \Rightarrow \ [\xi(\v_{_{2}})-\epsilon-t+1,1]\notin\v_{_{1}} \ \Rightarrow \ [0,\xi(\v_{_{2}})-\epsilon-t+1[\in\v_{_{1}}
\end{align*}
and, since $\xi(\v_{_{2}})-\epsilon<\xi(\v_{_{2}})$, then $[\xi(\v_{_{2}})-\epsilon,1]\in\v_{_{2}}$. Hence $[0,\xi(\v_{_{2}})-\epsilon-t+1[\times[\xi(\v_{_{2}})-\epsilon,1]\in\w$, but one can see that $[0,\xi(\v_{_{2}})-\epsilon-t+1[\times[\xi(\v_{_{2}})-\epsilon,1]\subseteq R_{_{t}}$, a contradiction.

Thus, for each $\w\in U([0,1]\times[0,1])$, $\xi\cdot U\vec{\circledast}(\w)=\xi(\v_{_{1}})\circledast\xi(\v_{_{2}})$ and the diagram is commutative. We have proved in the first item that the ultrafilter functor $U$ satisfies inequality (\ref{ineq_prod_alex}) for every set $X$. Therefore we conclude that Alexandroff $(\U,[0,1]_{_{\odot}})$-spaces form a cartesian closed subcategory of $\Cat{(\U,[0,1]_{_{\odot}})}$. Furthermore, Alexandroff $(\U,[0,1]_{_{\odot}})$-spaces are precisely the image by the functor $A^{\circ}\colon\categ{B}_{_{1}}\Met\to\Cat{(\U,[0,1]_{_{\odot}})}$ of the Alexandroff bounded-by-1 metric spaces. Explicitly, if $(X,d)$ is an Alexandroff bounded-by-1 metric space, then $(X,d_{_{\#}})$ is an Alexandroff $(\U,[0,1]_{_{\odot}})$-space, where 
\begin{align*}
    d_{_{\#}}(\x,x)=\inf_{A\in\x}\sup_{x'\in A}d(x',x);
\end{align*}
if $(X,a)$ is an Alexandroff $(\U,[0,1]_{_{\odot}})$-space, then $(X,a\cdot e_{_{X}})$ is an Alexandroff bounded-by-1 metric space, and, moreover, $a=d_{_{\#}}$, with $d=a\cdot e_{_{X}}$, so that, for each $\x\in UX$, $x\in X$,
\begin{align*}
    a(\x,x)=\inf_{A\in\x}\sup_{x'\in A}a(e_{_{X}}(x'),x).
\end{align*}

\noindent (3) For the category $\NA\text{-}\App$ we cannot apply Corollary \ref{sub_alex} nor Lemma 2, since the diagram
\begin{align*}\xymatrix{U([0,\infty]\times[0,\infty]) \ar[rr]^(0.55){U\vec{\ovee}} \ar[d]_{\left\langle\xi\cdot U\pi_{_{1}},\xi\cdot U\pi_{_{2}}\right\rangle} & & U[0,\infty] \ar[d]^{\xi} \\ [0,\infty]\times[0,\infty] \ar[rr]_(0.55){\vec{\ovee}} & & [0,\infty]}\end{align*}
is not commutative, where $\xi$ is defined as in the first item. To prove this fact, let $0<v<\infty$ and consider the filters on $[0,\infty]$ defined by
\begin{gather*}
    \mathfrak{f}_{_{1}}=\{A\subseteq[0,\infty] \ | \ \exists \ u<v; ]u,v]\subseteq A \ \text{or} \ [0,v]\subseteq A\} 
\shortintertext{and} 
    \mathfrak{f}_{_{2}}=\{A\subseteq[0,\infty] \ | \ \exists \ u>v, [0,u]\subseteq A \ \text{or} \ ]v,u]\subseteq A\}.
\end{gather*}
There exist ultrafilters $\mathcal{\v}_{_{1}}\supseteq\mathfrak{f}_{_{1}}$ and $\mathcal{\v}_{_{2}}\supseteq\mathfrak{f}_{_{2}}$; moreover, one can see that $\xi(\v_{_{1}})=\xi(\v_{_{2}})=v$, whence $\xi(\v_{_{1}})\ovee\xi(\v_{_{2}})=0$. Since the ultrafilter monad satisfies (BC) and $U\one=\one$, the diagram
\begin{align*}
    \xymatrix{U([0,\infty]\times[0,\infty]) \ar[d]_{U\pi_{_{1}}} \ar[r]^(0.62){U\pi_{_{2}}} & U[0,\infty] \ar[d] \\ U[0,\infty] \ar[r] & \one}
\end{align*}
is a weak pullback, hence there exists $\w\in U([0,\infty]\times[0,\infty])$ such that $U\pi_{_{i}}(\w)=\v_{_{i}}$, $i=1,2$. Now $\xi\cdot U\vec{\ovee}(\w)=\inf\{u\in[0,\infty] \ | \ D_{_{u}}\in\w\}$, where 
\begin{align*}
    D_{_{u}}=\{(p,q)\in[0,\infty]\times[0,\infty] \ | \ p\geq q \ \text{or} \ q\leq u\},
\end{align*}
which is depicted in the graphic below. \vspace{0.2cm}
\begin{center}
    \includegraphics[scale=0.5]{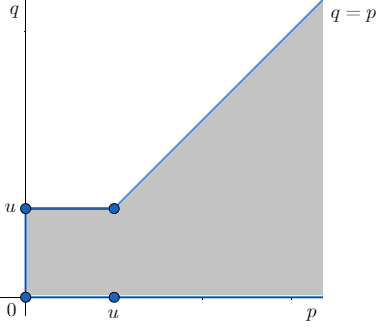}
\end{center} 
Consider an element $u<v$. Since $[0,v]\in\v_{_{1}}$ and $]v,\infty]\in\v_{_{2}}$, we have $[0,v]\times]v,\infty]\in\w$, and because $[0,v]\times]v,\infty]\subseteq\left([0,\infty]\times[0,\infty]\right)\backslash D_{_{u}}$, then $D_{_{u}}\notin\w$. Hence, if $D_{_{u}}\in\w$, then $v\leq u$, so that $\xi(\v_{_{1}})\ovee\xi(\v_{_{2}})=0<v\leq\xi\cdot U\vec{\ovee}(\w)$. 
\begin{center}
    \includegraphics[scale=0.5]{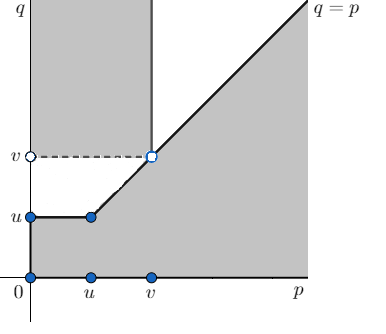}
\end{center}


\subsection{Exponentiably generated \texorpdfstring{$(\T,\V)$}{(T,V)}-spaces} 

Consider the class $\C$ of exponentiable spaces. One readily checks that property {\bf (EP)} of Subsection \ref{sub_cart_clos} holds. By Proposition \ref{sub_cgen} and Theorem \ref{sub_cart_clos}, we conclude that coequalizers of coproducts of exponentiable spaces form a cartesian closed subcategory of $\Cat{(\T,\V)}$. A criterion for exponentiability in $\Cat{(\T,\V)}$ is proved in \cite{MR3987969} (see (\ref{expo})).

In $\Ord$, $\C=\Ord_{_{\C}}=\Ord$. For $\Top$, exponentiable spaces are characterized as the {\it core-compact} spaces (see \cite{EH01} for a complete account on the matter). In \cite[Definition 1.8]{MR3227304} a space $(X,a)$ is said to be {\it core-compact} if $a\cdot Ta=a\cdot m_{_{X}}$, condition that implies $(X,a)$ to be $\otimes$-exponentiable \cite{MR2355608, MR3227304}, what in $\Top$ is equivalent to be (cartesian) exponentiable, since $\otimes=\wedge$ in the quantale $\two$. Exponentiably generated spaces in $\Top$, which are quotients of disjoint sums of core-compact spaces, are then called {\it core-compactly generated} \cite{MR2080286}. For an account on exponentiable metric spaces and exponentiable approach spaces see \cite{MR2259813} and \cite{MR3423697}, respectively.


\subsection{Injectively generated \texorpdfstring{$(\T,\V)$}{(T,V)}-spaces} 

Under the conditions of Theorem \ref{sub_inj_exp}, hence in particular for the categories of Table (\ref{tab1}), one can consider the class $\C$ of injective spaces. Then each element of $\C$ is exponentiable and binary products of injective spaces are again injective, hence condition {\bf (EP)} holds for $\C$. Proposition \ref{sub_cgen} and Theorem \ref{sub_cart_clos} imply that coequalizers of coproducts of injective spaces form a cartesian closed subcategory of $\Cat{(\T,\V)}$.

The injective spaces in the usual sense, that is, with ``$=$'' instead of ``$\simeq$'' in diagram (\ref{diag_ext}), also form a class satisfying condition {\bf (EP)}. In particular, in $\Top$ quotients of disjoint sums of retracts of powers $D^{I}$, $I$ a set and $D=(\{0,1,2\},\{\emptyset,\{0,1,2\},\{0,1\}\})$ \cite[Examples 9.3(4)]{MR1051419}, form a cartesian closed subcategory.


\section{Quasi-\texorpdfstring{$(\T,\V)$}{(T,V)}-spaces}\label{sec3}


\subsection{The category of quasi-spaces and quasi-continuous maps}\label{sub_qs}

Following the presentation of \cite{rsqsksd}, we generalize the work of \cite{MR0144300} from $\Top$ to $\Cat{(\T,\V)}$. Throughout the rest of the paper, let $\C$ denote the full subcategory of $\Cat{(\T,\V)}$ of compact Hausdorff spaces. We also assume the necessary conditions on $\T$ and $\V$ so that constant maps are continuous in $\Cat{(\T,\V)}$, and $\C$ is closed under finite coproducts, binary products, and equalizers. \vspace{0.2cm}

\noindent{\bf Definition} (1) For a map $\alpha\colon C\to X$ and a finite family $(\alpha_{_{i}})_{_{i\in I}}$ of maps $\alpha_{_{i}}\colon C_{_{i}}\to X$, with $C,C_{_{i}}\in\C$, one says that $\alpha$ {\it is covered by the family} $(\alpha_{_{i}})_{_{i\in I}}$ if there exists a surjective continuous map $\eta\colon\DS\coprod_{_{i\in I}}C_{_{i}}\to C$, where $\DS\coprod_{_{i\in I}}C_{_{i}}$ denotes the coproduct of the family $(C_{_{i}})_{_{i\in I}}$ in $\Cat{(\T,\V)}$, such that the triangle below is commutative.
\begin{align}\label{cover}\xymatrix@R=1em{\DS\coprod_{_{i\in I}}C_{_{i}} \ar[rrrrdd]^{\DS\coprod_{_{i\in I}}{\alpha_{_{i}}}} \ar[dd]_(0.6){\eta} & & & & \\ & & & & \\ C \ar[rrrr]_{\alpha} & & & & X}\end{align}
In particular, every map $\alpha$ is covered by itself. 

\noindent (2) A {\it quasi-$(\T,\V)$-space}, or simply a {\it quasi-space}, is a set $X$ together with, for each $C\in\C$, a set $Q(C,X)$ of functions from $C$ to $X$, whose elements are called {\it admissible maps}, satisfying the conditions: \\
(QS1) for all $C\in\C$, $Q(C,X)$ contains all constant maps, \\
(QS2) for all $C_{_{1}},C_{_{2}}\in\C$, if $h\colon C_{_{1}}\to C_{_{2}}$ is a $(\T,\V)$-continuous map and $\alpha\in Q(C_{_{2}},X)$, then $\alpha\cdot h\in Q(C_{_{1}},X)$, and \\
(QS3) for all $C\in\C$, $\alpha\in Q(C,X)$ if, and only if, $\alpha$ is covered by a family of admissible maps. \\ 
We denote a quasi-space by $(X,(Q(C,X))_{_{C\in\C}})$, or simply by $X$ when the corresponding {\it quasi-$(\T,\V)$-structure} $(Q(C,X))_{_{C\in\C}}$ is clear from the context. A map $f\colon X\to Y$ between quasi-spaces is said to be a {\it quasi-$(\T,\V)$-continuous} map, or simply a {\it quasi-continuous} map, if, for all $C\in\C$ and $\alpha\in Q(C,X)$, $f\cdot\alpha\in Q(C,Y)$; we denote the set of quasi-continuous maps from $X$ to $Y$ by $\Qs(X,Y)$. \vspace{0.2cm}

Identity maps and composition of quasi-continuous maps are quasi-continuous, so we have a category $\Qs\Cat{(\T,\V)}$. \vspace{0.2cm}

\noindent{\bf Remark} Although carrying a size illegitimacy {\bf --} its collection of objects do not form a class {\bf --} proved by Herrlich and Rajagopalan \cite{MR708599}, we still call $\Qs\Cat{(\T,\V)}$ a category. \vspace{0.2cm}

In the following items we discuss some properties of quasi-spaces. \vspace{0.2cm}

\noindent{\bf (I)} {\it Each space $(X,a)$ is a quasi-space}: for each $C\in\C$, define 
\begin{align*}Q_{_{a}}(C,X)=\{\alpha\colon(C,c)\to(X,a) \ | \ \alpha \text{ is continuous}\}.\end{align*} 
Thus $Q_{_{a}}(C,X)$ satisfies (QS1) and (QS2). Let us verify (QS3). If $\alpha\in Q_{_{a}}(C,X)$, then it is trivially covered by itself. Consider $\alpha\colon C\to X$ covered by a family of admissible maps $(\alpha_{_{i}})_{_{i\in I}}$ as in (\ref{cover}). Observe that $\alpha\cdot\eta=\DS\coprod_{_{i\in I}}{\alpha_{_{i}}}\in Q_{_{a}}(\DS\coprod_{_{i\in I}}C_{_{i}},X)$ and then, Axiom of Choice granted, we conclude 
\begin{align*}\begin{array}{rclll} c & = & c\cdot T\eta\cdot(T\eta)^{\circ} & & (\text{$T\eta$ is a surjective map}) \\[0.1cm]
& \leq & c\cdot T\eta\cdot b^{\circ}\cdot b\cdot(T\eta)^{\circ} & & (\text{$(\amalg_{_{i}}C_{_{i}},b)$ is compact}) \\[0.1cm]
& \leq & c\cdot c^{\circ}\cdot\eta\cdot b\cdot(T\eta)^{\circ} & & (\text{$\eta$ is $(\T,\V)$-continuous}) \\[0.1cm]
& \leq & \eta\cdot b\cdot(T\eta)^{\circ} & & (\text{$(C,c)$ is Hausdorff}) \\[0.1cm]
& \leq & \eta\cdot(\alpha\cdot\eta)^{\circ}\cdot a\cdot T(\alpha\cdot\eta)\cdot(T\eta)^{\circ} & & (\text{$\alpha\cdot\eta$ is $(\T,\V)$-continuous}) \\[0.1cm]
& = & \eta\cdot\eta^{\circ}\cdot\alpha^{\circ}\cdot a\cdot T\alpha\cdot T\eta\cdot(T\eta)^{\circ} & \\[0.1cm]
& \leq & \alpha^{\circ}\cdot a\cdot T\alpha & & (\text{$\eta$ and $T\eta$ are maps}).\end{array}\end{align*}
Observe that, by the same argument, one can conclude the following lemma (see also \cite[V-4.3.1]{MR3307673}). \vspace{0.2cm}

\noindent{\bf Lemma} {\it Every surjective $(\T,\V)$-continuous map from a compact $(\T,\V)$-space to a Hausdorff $(\T,\V)$-space is a proper $(\T,\V)$-continuous map, which is also a quotient map.} \vspace{0.2cm}

For details on proper $(\T,\V)$-continuous maps see \cite{MR2107395} and also \cite[V-3]{MR3307673}. The quasi-space $(X,(Q_{_{a}}(C,X))_{_{C\in\C}})$ is said to be {\it associated} with the space $(X,a)$. Moreover, each continuous map $f\colon(X,a)\to(Y,b)$ is quasi-continuous with respect to the associated quasi-spaces: for each $C\in\C$ and admissible $\alpha\colon C\to X$ (that is, $\alpha$ continuous), $f\cdot\alpha\colon C\to Y$ is continuous, hence admissible by definition. However, in general this inclusion of $\Cat{(\T,\V)}$ into $\Qs\Cat{(\T,\V)}$ is not full. \vspace{0.2cm}

\noindent{\bf (II)} {\it Different $(\T,\V)$-structures on the set $X$ may originate the same quasi-space}: take $(X,a)$ and the corresponding $\C$-generated space $(X,a^{c})$ defined in Section \ref{sec2}. By Lemma \ref{sub_cgen}, for each $C\in\C$, a map $\alpha\colon C\to(X,a)$ is continuous if, and only if, $\alpha\colon C\to(X,a^{c})$ is continuous. Actually, the $\C$-generated structure on $X$ is the least one inducing the same associated quasi-space, for if $\overline{a}$ is such a structure, then, in particular, each probe over $(X,a)$ is a continuous map $\alpha_{_{i}}\colon(X_{_{i}},a_{_{i}})\to(X,\overline{a})$, what implies that the identity $1_{_{X}}\colon(X,a^{c})\to(X,\overline{a})$ is continuous, hence $a^{c}\leq\overline{a}$.\vspace{0.2cm} 

\noindent{\bf (III)} {\it In general, there exist quasi-spaces in $\Qs\Cat{(\T,\V)}$ wich are not associated with any space in $\Cat{(\T,\V)}$}. For $\Top$, this is shown in \cite[Lemma 5.5]{MR0144300}. Let us study the case of $\Cat{\V}$. As seen in the Examples of Section \ref{sec1}, the compact Hausdorff $\V$-spaces coincide with the discrete $\V$-spaces, that is, spaces of the form $(C,1_{_{C}})$, $C\in\Set$. Moreover, for each $\V$-space $(X,a)$, any map $\alpha\colon(C,1_{_{C}})\to(X,a)$ is $\V$-continuous, whence the associated quasi-$\V$-structure on $X$ is given by, for each $C\in\Set$, 
\begin{align*}Q_{_{a}}(C,X)=\{\text{maps} \ \alpha\colon C\to X\}.\end{align*} 
Therefore, if $X$ is not a singleton, the quasi-$\V$-structure on $X$ defined by 
\begin{align*}Q'(C,X)=\{\text{constant maps} \ \alpha\colon C\to X\},\end{align*} 
for each $C\in\Set$, is not associated with any $\V$-space $(X,a)$ in $\Cat{\V}$.\vspace{0.2cm}

\noindent{\bf (IV)} {\it Admissible maps of the associated quasi-space, quasi-continuous maps, and $(\T,\V)$-continuous maps coincide when considering compact Hausdorff $(\T,\V)$-spaces.} Consider $(C,c)\in\C\subseteq\Cat{(\T,\V)}$ and its associated quasi-structure $(Q_{_{c}}(B,C))_{_{B\in\C}}=(\C(B,C))_{_{B\in\C}}$. Let $(X,(Q(B,X))_{_{B\in\C}})$ be a quasi-space, and $\alpha\colon C\to X$ be a map. If $\alpha$ is quasi-continuous, then $\alpha\cdot1_{_{C}}=\alpha\in Q(C,X)$, since $1_{_{C}}\in Q_{_{c}}(C,C)$. Conversely, if $\alpha$ is an admissible map from $C$ to $X$, then, for each $B\in\C$ and $\beta\in Q_{_{c}}(B,C)=\C(B,C)$, $\alpha\cdot\beta\in Q(B,X)$, by (QS2), so $\alpha$ is quasi-continuous. Hence 
\begin{align*}Q(C,X)=\Qs(C,X).\end{align*} 
In particular, if $X=(D,d)\in\C$, quasi-continuous maps between the associated quasi-spaces $C$ and $D$ coincide with the admissible maps $Q_{_{d}}(C,D)$, which are all continuous maps from $C$ to $D$ by definition, that is, 
\begin{align*}Q_{_{d}}(C,D)=\Qs(C,D)=\Cat{(\T,\V)}(C,D).\end{align*}

\noindent{\bf (V)} {\it Discrete and indiscrete quasi-structures}. For each set $X$ we have a {\it discrete} quasi-structure given by, for each $C\in\C$, 
\begin{align*}Q_{_{{\rm dis}}}(C,X)=\{\text{constant maps from $C$ to $X$}\}.\end{align*} 
Denoting $(X,(Q_{_{{\rm dis}}}(C,X))_{_{C\in\C}})$ by $DX$, we see that, for each quasi-space $(Y,(Q(C,Y))_{_{C\in\C}})$, each map $f\colon X\to Y$ is a quasi-continuous map $f\colon DX\to(Y,(Q(C,Y))_{_{C\in\C}})$. Analogously, we can endow the set $X$ with an {\it indiscrete} quasi-structure: for each $C\in\C$,  
\begin{align*}Q_{_{{\rm ind}}}(C,X)=\Set(C,X).\end{align*} 
Denoting $(X,(Q_{_{{\rm ind}}}(C,X))_{_{C\in\C}})$ by $IX$, for each quasi-space $(Y,(Q(C,Y))_{_{C\in\C}})$, each map $f\colon Y\to X$ is a quasi-continuous map $f\colon(Y,(Q(C,Y))_{_{C\in\C}})\to IX$. We observe that the indiscrete quasi-space is associated with the indiscrete space $(X,\top)$, because, for each $C\in\C$, every map $f\colon C\to(X,\top)$ is continuous. The same does not happen in general with the discrete quasi-space, as one verifies, for example in $\Top$, that there exist continuous maps from compact Hausdorff spaces to discrete spaces that are not constant. \vspace{0.2cm}

With the definitions in item {\bf (V)} we prove the following: \vspace{0.2cm}

\noindent{\bf Proposition} {\it The forgetful functor $|\text{-}|\colon\Qs\Cat{(\T,\V)}\to\Set$ has left and right adjoints, and it is represented by the singleton quasi-space.}

\begin{proof}{\rm For the left adjoint, define $D\colon\Set\to\Qs\Cat{(\T,\V)}$ assigning to each set $X$ the discrete quasi-space $DX$, and on maps $D$ is the identity. As observed before, for each quasi-space $(Y,(Q(C,Y))_{_{C\in\C}})$, each $\Set$-map $f\colon X\to Y$ is a quasi-continuous map $f\colon DX\to(Y,(Q(C,Y))_{_{C\in\C}})$, whence 
\begin{align*}\Qs(DX,(Y,(Q(C,Y))_{_{C\in\C}}))=\Set(X,Y)=\Set(X,|(Y,(Q(C,Y))_{_{C\in\C}})|).\end{align*}

For the right adjoint, define $I\colon\Set\to\Qs\Cat{(\T,\V)}$ assigning to each set $X$ the indiscrete quasi-space $IX$, and leaving the morphisms unchanged; for each quasi-space $(Y,(Q(C,Y))_{_{C\in\C}})$, a map $f\colon X\to Y$ is a quasi-continuous map $f\colon(X,(Q(C,X))_{_{C\in\C}})\to IY$, so 
\begin{align*}\Set(|(X,(Q(C,X))_{_{C\in\C}})|,Y)=\Set(X,Y)=\Qs((X,(Q(C,X))_{_{C\in\C}}),IY).\end{align*}

Since $|\text{-}|$ has a left adjoint $D$, it is represented by $D\one$, which coincides with $I\one$, and is given by the singleton $\one=\{*\}$ endowed with the quasi-structure defined by 
\begin{align*}Q(C,\one)=\{!_{_{C}}\colon C\to\one\},\end{align*} 
for each $C\in\C$.}\end{proof}


\subsection{\texorpdfstring{$\Qs\Cat{(\T,\V)}$}{Qs(T,V)-Cat} is topological over \texorpdfstring{$\Set$}{Set}}\label{sub_lim_colim}

Given a quasi-space $(X,(Q(C,X))_{_{C\in\C}})$ and a subset $A\subseteq X$, we can consider the {\it subspace quasi-structure} on $A$, which is given by, for each $C\in\C$,  
\begin{align}\label{sub_quasi_struc}\alpha\in Q(C,A) \ \Longleftrightarrow \ i_{_{A}}\cdot\alpha\in Q(C,X),\end{align} 
where $i_{_{A}}\colon A\hookrightarrow X$ is the inclusion map. When $A$ is endowed with this structure, $i_{_{A}}$ becomes a quasi-continuous map which is also $|\text{-}|$-initial, with $|\text{-}|\colon\Qs\Cat{(\T,\V)}\to\Set$ the forgetful functor.

Furthermore, for a quasi-space $(X,(Q(C,X))_{_{C\in\C}})$ and a surjective map $f\colon X\to Y$, we can define a {\it quotient} quasi-structure by: for each $C\in\C$, $\alpha\in Q(C,Y)$ if there exist a surjective map $f'\colon C'\to C$ in $\C$, and a map $\alpha'\in Q(C',X)$ such that the square below is commutative. 
\begin{align}\label{quo}\xymatrix{C' \ar[r]^{f'} \ar[d]_{\alpha'} & C \ar[d]^{\alpha} \\ X \ar[r]_{f} & Y}\end{align}

One can check that the latter structure satisfies (QS1) and (QS3). To verify (QS2), take \linebreak $\alpha\in Q(C,Y)$ and $h\colon B\to C$ be a continuous map, with $B,C\in\C$. By definition, there exist a surjective map $f'\colon C'\to C$ and a map $\alpha'\in Q(C',X)$, for $C'\in\C$, as in (\ref{quo}). Take the pullback of $f'$ along $h$ as in the diagram below.
\begin{align*}\xymatrix{B\times_{_{C}}C' \ar[rr]^{\pi_{_{1}}} \ar[d]_{\pi_{_{2}}} \ar@<-0.0ex>@{}[dr]|(0.3){\text{\pigpenfont J}} & & B \ar[d]^{h} \\ C' \ar[rr]_{f'} \ar[d]_{\alpha'} & & C \ar[d]^{\alpha} \\ X \ar[rr]_{f} & & Y}\end{align*}
Since $f'$ is a surjective $\Set$-map, $\pi_{_{1}}$ is surjective too. Also, by our assumptions, $B\times_{_{C}}C'\in\C$, and because $\alpha'\in Q(C',X)$, we have $\alpha'\cdot\pi_{_{2}}\in Q(B\times_{_{C}}C',X)$. 

When $Y$ is endowed with the quotient quasi-structure with respect to the surjection $f\colon X\to Y$, the map $f$ becomes not only quasi-continuous, but also a $|\text{-}|$-final morphism: if $g\colon Y\to Z$ is a map such that $g\cdot f\colon X\to Z\in\Qs(X,Z)$, for $(Z,(Q(C,Z))_{_{C\in\C}})\in\Qs\Cat{(\T,\V)}$, then, for each $C\in\C$ and $\alpha\in Q(C,Y)$, there exist a surjection $f'\colon C'\to C$, and a map $\alpha'\in Q(C',X)$ commuting the square in (\ref{quo}), hence 
\begin{align*}g\cdot\alpha\cdot f'=g\cdot f\cdot\alpha'\in Q(C',Z),\end{align*}
so the map $g\cdot\alpha$ is covered by an admissible map:
\begin{align*}\xymatrix{C' \ar[rrrd]^{g\cdot f\cdot\alpha'} \ar[d]_{f'} & & & \\ C \ar[rrr]_{g\cdot\alpha} & & & Z,}\end{align*}
whence $g\cdot\alpha\in Q(C,Z)$, and so $g\in\Qs(Y,Z)$.

The constructions above lead us to the following result. \vspace{0.2cm}

\noindent{\bf Proposition} {\it The forgetful functor $|\text{-}|\colon\Qs\Cat{(\T,\V)}\to\Set$ is topological.} \vspace{0.2cm}

\noindent{\it Proof.} Let $(f_{_{j}}\colon X\to|(X_{_{j}},(Q(C,X_{_{j}}))_{_{C\in\C}})|)_{_{j\in J}}$ be a source in $\Set$. To construct a $|\text{-}|$-initial lifting
\begin{align*}(f_{_{j}}\colon(X,(Q(C,X))_{_{C\in\C}})\to(X_{_{j}},(Q(C,X_{_{j}}))_{_{C\in\C}}))_{_{j\in J}}\end{align*} 
for this source, define, for each $C\in\C$, 
\begin{align*}\alpha\in Q(C,X) \ \Longleftrightarrow \ \forall j\in J, \ f_{_{j}}\cdot\alpha\in Q(C,X_{_{j}}).\end{align*} 
Properties (QS1) and (QS2) are immediately satisfied. To check (QS3), observe that if a map $\alpha\colon C\to X$ is covered by admissible maps $\alpha_{_{i}}\colon C_{_{i}}\to X$, $i\in I$ finite, then, for each $j\in J$, $f_{_{j}}\cdot\alpha$ is covered by the family of maps ($\beta_{_{ji}}=f_{_{j}}\cdot\alpha_{_{i}})_{_{i\in I}}$, which are admissible by definition of $Q(C_{_{i}},X)$.
\begin{align*}\xymatrix@=1.2em{\DS\coprod_{_{i\in I}}C_{_{i}} \ar[rrrrdd]^(0.55){\DS\coprod_{_{i\in I}}{\alpha_{_{i}}}} \ar[dd]_{\eta} \ar@/^2.5pc/[rrrrrrdd]^(0.8){\DS\coprod_{_{i\in I}}{\beta_{_{ji}}}} & & & & & & \\ & & & & & & \\ C \ar[rrrr]_{\alpha} & & & & X \ar[rr]_{f_{_{j}}} & & X_{_{j}}}\end{align*}
Hence, for all $j\in J$, $f_{_{j}}\cdot\alpha\in Q(C,X_{_{j}})$, so $\alpha\in Q(C,X)$.  

It is straightforward to verify that the latter lifting is indeed $|\text{-}|$-initial and uniqueness follows from amnesticity of $|\text{-}|$ \cite[Proposition 21.5]{MR1051419}. \qed \vspace{0.2cm}

Among other properties, the latter result implies completeness and cocompleteness of $\Qs\Cat{(\T,\V)}$ \cite[Proposition 13.15]{MR1051419}; following this reference we describe limits and colimits. \vspace{0.2cm}

\noindent\textbf{\textit{Limits.}} Let $\A$ be a small category and $\mD\colon\A\to\Qs\Cat{(\T,\V)}$ be a diagram. First construct the limit in $\Set$ of $|\text{-}|\cdot\mD\colon\A\to\Set$, that we denote by $(p_{_{A}}\colon X\to|\mD A|)_{_{A\in{\rm Obj}\A}}$, and then take the $|\text{-}|$-initial lifting of this source, described in the previous Proposition. In particular, the product of a family $((X_{_{i}},(Q(C,X_{_{i}}))_{_{C\in\C}}))_{_{i\in I}}$ of quasi-spaces is given by the set $\DS\prod_{_{i\in I}}{X_{_{i}}}$ endowed with the quasi-structure: for each $C\in\C$, 
\begin{align*}\alpha\in Q(C,\DS\prod_{_{i\in I}}{X_{_{i}}}) \ \Longleftrightarrow \ \forall i\in I, \ \pi_{_{i}}\cdot\alpha\in Q(C,X_{_{i}}),\end{align*} 
where the $\pi_{_{i}}$'s are product projections. One can see that for the empty family, the product is given by the singleton $\one$ endowed with the quasi-structure: for each $C\in\C$, $Q(C,\one)=\{!_{_{C}}\colon C\to\one\}$, which was described in Proposition \ref{sub_qs}. As for equalizers of quasi-continuous maps $f,g\colon X\to Y$, endow the set $E=\{x\in X \ | \ f(x)=g(x)\}\subseteq X$ with the subspace quasi-structure.\vspace{0.2cm}  

\noindent\textbf{\textit{Colimits.}} For a diagram $\mD\colon\A\to\Qs\Cat{(\T,\V)}$, we form the colimit in $\Set$ of $|\text{-}|\cdot\mD\colon\A\to\Set$, denoted by $(j_{_{A}}\colon|\mD A|\to X)_{_{A\in{\rm Obj}\A}}$, and then we take the $|\text{-}|$-final lifting of this sink. The quasi-structure on $X$ is given by: for each $C\in\C$, $\alpha\in Q(C,X)$ if, and only if, $\alpha$ is covered by a family $(\alpha_{_{i}})_{_{i\in I}}$ such that each $\alpha_{_{i}}$ factorizes through a colimit inclusion $j_{_{A_{_{i}}}}\colon\mD A_{_{i}}\to X$ and an admissible map $\beta_{_{i}}\in Q(C_{_{i}},\mD A_{_{i}})$.
\begin{align*}\xymatrix@=1em{\coprod{C_{_{i}}} \ar[rrrdd]^{\coprod{\alpha_{_{i}}}} \ar[dd]_{\eta} & & & C_{_{i}} \ar@{_{(}->}[lll]^{} \ar[dd]^{\alpha_{_{i}}} \ar[rrd]^{\beta_{_{i}}} & & \\ & & & & & \mD A_{_{i}} \ar[lld]^{j_{_{A_{_{i}}}}} \\ C \ar[rrr]_{\alpha} & & & X & & }\end{align*}
In particular cases, we can reduce this quasi-structure (see \cite{rsqsksd}). The coproduct of a family $((X_{_{i}},(Q(C,X_{_{i}}))_{_{C\in\C}}))_{_{i\in I}}$ of quasi-spaces is given by the disjoint union $\dot{\bigcup}{X_{_{i}}}$ endowed with the quasi-structure: for $C\in\C$, $\alpha\in Q(C,\dot{\bigcup}{X_{_{i}}})$ if, and only if, $\alpha$ is covered by a family $(j_{_{i_{_{k}}}}\cdot\beta_{_{k}})_{_{k\in K}}$, $K$ a finite set, with $\eta=1_{_{C}}$, $j_{_{i_{_{k}}}}\colon X_{_{i_{_{k}}}}\to\dot{\bigcup}{X_{_{i}}}$ the coproduct inclusion, and $\beta_{_{k}}\in Q(C_{_{k}},X_{_{i_{_{k}}}})$.
\begin{align*}\xymatrix@=1em{\coprod{C_{_{k}}} \ar[rrrdd]^{\coprod{(j_{_{i_{_{k}}}}\cdot\beta_{_{k}})}} \ar[dd]_{1_{_{C}}} & & & C_{_{k}} \ar@{_{(}->}[lll]^{} \ar[dd]^{} \ar[rrd]^{\beta_{_{k}}} & & \\ & & & & & X_{_{i_{_{k}}}} \ar[lld]^{j_{_{i_{_{k}}}}} \\ C \ar[rrr]_{\alpha} & & & \dot{\bigcup}{X_{_{i}}} & & }\end{align*} 
The initial object is then given by $\emptyset$ endowed with the quasi-structure: 
\begin{align*}Q(C,\emptyset)=\left\{\begin{array}{l} \emptyset, \ \text{if} \ C\neq\emptyset \\ \{1_{_{\emptyset}}\}, \ \text{otherwise.}\end{array}\right.\end{align*} 
As for coequalizers of quasi-continuous maps $f,g\colon X\to Y$, consider in $Y$ the smallest equivalence relation containing the pairs $(f(x),g(x))$, for $x\in X$, and endow $\tilde{Y}=Y/\hspace{-0.1cm}\sim$ with the quotient quasi-structure with respect to the projection map $p_{_{Y}}\colon Y\to\tilde{Y}$. 


\subsection{\texorpdfstring{$\Qs\Cat{(\T,\V)}$}{Qs(T,V)-Cat} is cartesian closed} 

In general $(\T,\V)$-spaces and $(\T,\V)$-continuous maps do not form a cartesian closed category. Hence this property is desirable for a supercategory of $\Cat{(\T,\V)}$.

The natural candidate for an exponential of quasi-spaces $X$ and $Y$ is $\Qs(X,Y)$. Consider the evalutation map ${\rm ev}\colon\Qs(X,Y)\times X\to Y$, $(f,x)\mapsto f(x)$. First we wish to define a quasi-structure on $\Qs(X,Y)$ such that ${\rm ev}$ is a quasi-continuous map, that is, for each $\gamma\in Q(C,\Qs(X,Y)\times X)$, with $C\in\C$, ${\rm ev}\cdot\gamma\in Q(C,Y)$. Hence, for each $\beta\in Q(C,\Qs(X,Y))$ and $\alpha\in Q(C,X)$, the composite ${\rm ev}\cdot\left\langle \beta,\alpha\right\rangle$ must belong to $Q(C,Y)$. 

Under this intuition and keeping in mind conditions (QS1), (QS2), and (QS3), define, for each $C\in\C$, $\beta\in Q(C,\Qs(X,Y))$ if for each $(\T,\V)$-continuous map $h\colon B\to C$, for $B\in\C$, and each $\alpha\in Q(B,X)$, the map ${\rm ev}\cdot\left\langle\beta\cdot h,\alpha\right\rangle\colon B\to Y$ belongs to $Q(B,Y)$. 

The latter data indeed define a quasi-structure and to verify, for instance, (QS3), take a map $\beta\colon C\to\Qs(X,Y)$ covered by a family of admissible maps $(\beta_{_{i}})_{_{i\in I}}$ as in (\ref{cover}). For a continuous map $h\colon B\to C$, with $B\in\C$, and $\alpha\in Q(B,X)$, form the following pullbacks
\begin{align*}\xymatrix@=1em{(\coprod_{_{i}}{C_{_{i}}})\times_{_{C}}B \ar[rrr]^(0.6){\pi_{_{B}}} \ar[dd]_{\pi_{_{\coprod_{_{i}}{C_{_{i}}}}}} \ar@<-0.0ex>@{}[dr]|(0.4){\text{\pigpenfont J}} & & & B \ar[dd]^{h} \\ & & & \\ \coprod_{_{i}}{C_{_{i}}} \ar[rrr]_(0.6){\eta} & & & C}\hspace{2cm}\xymatrix@=1em{C_{_{i}}\times_{_{C}}B \ar[rrr]^(0.6){\pi^{i}_{_{B}}} \ar[dd]_{\pi_{_{C_{_{i}}}}} \ar@<-0.0ex>@{}[dr]|(0.4){\text{\pigpenfont J}} & & & B \ar[dd]^{h} \\ & & & \\ C_{_{i}} \ar[rrr]_(0.6){\eta_{_{i}}} & & & C,}\end{align*}
where, for each $i\in I$, $\eta_{_{i}}=\eta\cdot j_{_{i}}$, with $j_{_{i}}\colon C_{_{i}}\hookrightarrow\coprod_{_{i}}{C_{_{i}}}$ the coproduct inclusion, and consider the commutative diagram
\begin{align*}\xymatrix@=1.2em{\coprod_{_{i}}{(C_{_{i}}\times_{_{C}}B)} \ar[d]_{\mu} \ar@{^{(}->}[rr] \ar@/^4.5pc/[rrrrrrdd]^(0.8){\coprod_{_{i}}{\gamma_{_{i}}}} & & \coprod_{_{i}}{(C_{_{i}}\times B)} \ar[d]_{\cong} \ar@<1ex>@/^1.5pc/[rrdd]^(0.7){\coprod_{_{i}}{(\beta_{_{i}}\times\alpha)}} & & & & \\ (\coprod_{_{i}}{C_{_{i}}})\times_{_{C}}B \ar[dd]_{\pi_{_{B}}} \ar@{^{(}->}[rr] & & (\coprod_{_{i}}{C_{_{i}}})\times B \ar[d]_{\eta\times1_{_{B}}} \ar[rrd]^{(\coprod_{_{i}}{\beta_{_{i}}})\times\alpha} & & & & \\ & & C\times B \ar[rr]_(0.4){\beta\times\alpha} & & \Qs(X,Y)\times X \ar[rr]_(0.6){{\rm ev}} & & Y, \\ B \ar@/_2pc/[rrrrrru]_{{\rm ev}\cdot\left\langle\beta\cdot h,\alpha\right\rangle} &  & & & & & }\end{align*} 
where $\mu$ is the surjective $(\T,\V)$-continuous map
\begin{align*}\begin{array}{rcl}\mu\colon\coprod_{_{i}}{(C_{_{i}}\times_{_{C}}B)} & \longrightarrow & (\coprod_{_{i}}{C_{_{i}}})\times_{_{C}}B \\ ((c_{_{i}},b),i) & \longmapsto & ((c_{_{i}},i),b).\end{array}\end{align*} 
We observe that we also use distributivity of $\Cat{(\T,\V)}$. The map ${\rm ev}\cdot\left\langle\beta\cdot h,\alpha\right\rangle$ is then covered by the family of maps $(\gamma_{_{i}})_{_{i\in I}}$, where, for each $(c_{_{i}},b)\in C_{_{i}}\times_{_{C}}B$, 
\begin{align*}\gamma_{_{i}}(c_{_{i}},b)=\beta_{_{i}}(\eta_{_{i}}(c_{_{i}}))(\alpha(b))=\beta_{_{i}}(\eta_{_{i}}\cdot\pi_{_{C_{_{i}}}}(c_{_{i}},b))(\alpha\cdot\pi^{i}_{_{B}}(c_{_{i}},b)),\end{align*} 
so $\gamma_{_{i}}={\rm ev}\cdot\left\langle\beta_{_{i}}\cdot\eta_{_{i}}\cdot\pi_{_{C_{_{i}}}},\alpha\cdot\pi^{i}_{_{B}}\right\rangle$ is admissible, since $\beta_{_{i}}$ is admissible, $\eta_{_{i}}\cdot\pi_{_{C_{_{i}}}}\colon C_{_{i}}\times_{_{C}}B\to C$ is continuous, and $\alpha$ admissible implies $\alpha\cdot\pi^{i}_{_{B}}$ admissible.

Choosing $h=1_{_{C}}$ in the definition of the quasi-structure, the map ${\rm ev}\colon\Qs(X,Y)\times X\to Y$ proves to be quasi-continuous. Furthermore, for each quasi-continuous map $f\colon Z\times X\to Y$, for $Z\in\Qs\Cat{(\T,\V)}$, there exists a unique $\Set$-map $\overline{f}\colon Z\to\Qs(X,Y)$, the transpose of $f$, such that ${\rm ev}\cdot(\overline{f}\times1_{_{X}})=f$. We verify next that $\overline{f}$ is quasi-continuous. 
\begin{align*}\xymatrix@R=1em{\Qs(X,Y) & \Qs(X,Y)\times X \ar[rr]^(0.6){{\rm ev}} & & Y \\ & & & \\ Z \ar@{.>}[uu]^{\exists \ ! \ \overline{f}} & Z\times X \ar[uu]^{\overline{f}\times1_{_{X}}} \ar[uurr]_(0.6){f} & & }\end{align*}
Let $\gamma\in Q(C,Z)$, with $C\in\C$; we wish to prove that $\overline{f}\cdot\gamma\in Q(C,\Qs(X,Y))$. For that, let $h\colon B\to C$ be a continuous map, for $B\in\C$, and $\alpha\in Q(B,X)$. Then $\gamma\cdot h\in Q(B,Z)$ and $\left\langle \gamma\cdot h,\alpha\right\rangle\in Q(B,Z\times X)$, whence $f\cdot\left\langle \gamma\cdot h,\alpha\right\rangle\in Q(B,Y)$. The result follows from the equalities: for each $b\in B$, 
\begin{align*}\begin{array}{rcl} {\rm ev}\cdot\left\langle\overline{f}\cdot\gamma\cdot h,\alpha\right\rangle(b) & = & \overline{f}\cdot\gamma\cdot h(b)(\alpha(b)) \\[0.1cm]
& = & \overline{f}(\gamma\cdot h(b))(\alpha(b)) \\[0.1cm]
& = & f(\gamma\cdot h(b),\alpha(b)) \\[0.1cm]
& = & f\cdot\left\langle \gamma\cdot h,\alpha\right\rangle(b).\end{array}\end{align*}
We have proved the following: \vspace{0.2cm}

\noindent{\bf Theorem} {\it $\Qs\Cat{(\T,\V)}$ is cartesian closed.} \vspace{0.2cm}

\noindent{\bf Examples} (1) Let us begin with $\Cat{\V}$. $C$ is a compact Hausdorff $\V$-space if, and only if, $C$ is discrete, that is, of the form $(C,1_{_{C}})$; this way, a quasi-$\V$-space is a set $X$ with, for each $(C,1_{_{C}})\in\C$, a set of maps $Q(C,X)$ satisfying conditions (QS1), (QS2), and (QS3). Moreover, a quasi-space $(X,(Q_{_{a}}(C,X))_{_{C\in\C}})$ associated with a $\V$-space $(X,a)$ necessarily satisfies 
\begin{align*}Q_{_{a}}(C,X)=\{\text{maps} \ \alpha\colon C\to X\}=\Set(C,X).\end{align*} 
Therefore, associated quasi-spaces coincide with indiscrete quasi-spaces (items {\bf (III)} and {\bf (V)} of Subsection \ref{sub_qs}). In particular, this happens for the quantales $\two$, $\categ{P}_{_{+}}$, $\categ{P}_{_{{\rm max}}}$, and $[0,1]_{_{\odot}}$. The respective categories of quasi-spaces coincide: 
\begin{align*}\Qs\Ord=\Qs\Met=\Qs\UltMet=\Qs\categ{B}_{_{1}}\Met.\end{align*} 

\noindent(2) For $\Top$ we recover Spanier's category of quasi-topological spaces \cite{MR0144300}. For $\App$, $\NA$-$\App$, and $\Cat{(\U,[0,1]_{_{\odot}})}$, similar to the previous item, our observation in Section \ref{sec1} that
\begin{align*}\Cat{(\U,[0,1]_{_{\odot}})}_{_{{\rm CompHaus}}}\cong\NA\text{-}\App_{_{{\rm CompHaus}}}\cong\App_{_{{\rm CompHaus}}}\cong\Top_{_{{\rm CompHaus}}}\cong\Set^{\U}\end{align*}
allow us to conclude that
\begin{align*}\Qs\Cat{(\U,[0,1]_{_{\odot}})}=\Qs\NA\text{-}\App=\Qs\App=\Qs\Top.\end{align*} 


\section{Relationship between quasi-spaces and compactly generated spaces}

This relationship is studied for $\Top$ in \cite{rsqsksd}. We have commented in item {\bf (I)} of Subsection \ref{sub_qs} that the inclusion
\begin{align*}\Cat{(\T,\V)}\hookrightarrow\Qs\Cat{(\T,\V)}\end{align*} 
in general is not full. However, the situation is different when restricting ourselves to compactly generated spaces. 

Let $(X,a)$ be a compactly generated space; by definition, a map $f\colon(X,a)\to(Y,b)$, for $(Y,b)$ in $\Cat{(\T,\V)}$, is continuous if, and only if, for all continuous maps $\alpha\colon C\to Y$, with $C\in\C$, $f\cdot\alpha\colon C\to Y$ is continuous. Considering the quasi-spaces associated with $(X,a)$ and $(Y,b)$, then $f\colon X\to Y$ is continuous if, and only if, $f\colon X\to Y$ is quasi-continuous, that is,
\begin{align}\label{eq31}\Cat{(\T,\V)}(X,Y)=\Qs(X,Y).\end{align}
Conversely, if $(X,a)$ is a $(\T,\V)$-space such that, for each $(\T,\V)$-space $(Y,b)$, its associated quasi-space satisfies (\ref{eq31}), then, by the same reasoning, $(X,a)$ is compactly generated. 

We identify the elements of $\Cat{(\T,\V)}_{_{\C}}$ with their associated quasi-spaces, and denote the resulting subcategory of $\Qs\Cat{(\T,\V)}$ by $\C$-$\Cat{(\T,\V)}$. \vspace{0.2cm}

\noindent{\bf Proposition} {\it $\C$-$\Cat{(\T,\V)}$ is fully reflective in $\Qs\Cat{(\T,\V)}$.}\vspace{0.2cm}

\noindent{\it Proof.} For $X,Y\in\C$-$\Cat{(\T,\V)}$, by (\ref{eq31}), 
\begin{align*}\C\text{-}\Cat{(\T,\V)}(X,Y)=\Cat{(\T,\V)}(X,Y)=\Qs(X,Y),\end{align*} 
hence it is a full subcategory.

For each $(X,(Q(C,X))_{_{C\in\C}})$ in $\Qs\Cat{(\T,\V)}$, define $(X,a_{_{Q}})$, with $a_{_{Q}}$ the $(\T,\V)$-structure from the $|\text{-}|$-final lifting in $\Cat{(\T,\V)}$ of the sink $(\alpha\colon(C,c)\to X)_{_{C\in\C,\alpha\in Q(C,X)}}$. Consider the quasi-space $(X,(Q_{_{a_{_{Q}}}}(C,X))_{_{C\in\C}})$ associated with $(X,a_{_{Q}})$. 

Let $(Y,b)$ be a $(\T,\V)$-space and $f\colon X\to Y$ be a map. If $f\colon(X,a_{_{Q}})\to(Y,b)$ is continuous, then 
\begin{align*}f\colon(X,(Q_{_{a_{_{Q}}}}(C,X))_{_{C\in\C}})\to(Y,(Q_{_{b}}(C,Y))_{_{C\in\C}})\end{align*} 
is quasi-continuous. Conversely, if $f$ is quasi-continuous with respect to the associated quasi-structures, then, for each $C\in\C$ and each continuous map $\alpha\colon C\to(X,a_{_{Q}})$, the composite $f\cdot\alpha\colon C\to(Y,b)$ is continuous. Hence $f\colon(X,a_{_{Q}})\to(Y,b)$ is continuous, by definition of $a_{_{Q}}$. Therefore, for every $(\T,\V)$-space $(Y,b)$, $(X,a_{_{Q}})$ satisfies (\ref{eq31}), that is, $(X,(Q_{_{a_{_{Q}}}}(C,X))_{_{C\in\C}})$ belongs to $\C$-$\Cat{(\T,\V)}$.    

The identity map $1_{_{X}}\colon(X,(Q(C,X))_{_{C\in\C}})\to(X,(Q_{_{a_{_{Q}}}}(C,X))_{_{C\in\C}})$ is quasi-continuous, since each $\alpha\in Q(C,X)$, for $C\in\C$, is a continuous map $\alpha\colon C\to(X,a_{_{Q}})$. For each $(Y,(Q_{_{b}}(C,Y))_{_{C\in\C}})$ in $\C$-$\Cat{(\T,\V)}$, each quasi-continuous map $f\colon(X,(Q(C,X))_{_{C\in\C}})\to(Y,(Q_{_{b}}(C,Y))_{_{C\in\C}})$ induces a continuous $f\colon(X,a_{_{Q}})\to(Y,b)$, hence a quasi-continuous $f\colon(X,(Q_{_{a_{_{Q}}}}(C,X))_{_{C\in\C}})\to(Y,(Q_{_{b}}(C,Y))_{_{C\in\C}})$.
\begin{align*}\xymatrix@R=1em{(X,Q(C,X)) \ar[rr]^{1_{_{X}}} \ar[rrdd]_(0.45){f} & & (X,Q_{_{a_{_{Q}}}}(C,X)) \ar[dd]^{f} \\ & & \\ & & (Y,Q_{_{b}}(C,Y))}\end{align*}\qed\vspace{0.2cm}

We can draw the following explanatory diagram, where $\C$ stands for the class of compact Hausdorff $(\T,\V)$-spaces, as fixed in the beginning of Subsection \ref{sub_qs}.
\begin{align*}\xymatrix@R=1em{\Cat{(\T,\V)}_{_{\C}} \ar@{^{(}->}[rr]^{\top} \ar[dd]_{\cong} & & \Cat{(\T,\V)} \ar[dd]^{\text{non-full}} \ar@/_1.5pc/[ll] \\ & & \\ \C\text{-}\Cat{(\T,\V)} \ar@{^{(}->}[rr]_{\top} & & \Qs\Cat{(\T,\V)} \ar@/^1.5pc/[ll]}\end{align*}\vspace{0.1cm}

\section*{Acknowledgements}

This work was developed as part of the author's PhD thesis, under the supervision of Maria Manuel Clementino, whom the author thanks for proposing the problem of investigation, advising the whole study, and for several improvements made to the article.

\end{document}